\documentclass[12pt]{amsart}

\usepackage{latexsym, amssymb, amscd, amsfonts}

\usepackage{mathpple}
\usepackage{pstricks}
\usepackage{pst-poly}
\usepackage{multido}
\usepackage{pst-node}
\usepackage[all]{xypic}
\usepackage{ifthen}

\newboolean{showpicture}
\setboolean{showpicture}{false}

\newboolean{bourbaki}
\setboolean{bourbaki}{true}

\newboolean{draft}
\setboolean{draft}{true}

\def\gb{\mathfrak{b}}

\def\gg{\mathfrak{g}}
\def\gh{\mathfrak{h}}

\def\gl{\mathfrak{l}}

\def\gp{\mathfrak{p}}

\def\gs{\mathfrak{s}}

\def\gB{\mathfrak{B}}

\def\gG{\mathfrak{G}}
\def\gH{\mathfrak{H}}

\def\gK{\mathfrak{K}}
\def\gL{\mathfrak{L}}

\def\gN{\mathfrak{N}}

\def\gP{\mathfrak{P}}

\def\gZ{\mathfrak{Z}}

\def\C{\mathbb{C}}

\def\F{\mathbb{F}}

\def\Q{\mathbb{Q}}
\def\R{\mathbb{R}}

\def\Z{\mathbb{Z}}

\def\cA{\mathcal{A}}

\def\cC{\mathcal{C}}
\def\cD{\mathcal{D}}

\def\cG{\mathcal{G}}
\def\cH{\mathcal{H}}

\def\cL{\mathcal{L}}

\def\cO{\mathcal{O}}
\def\cP{\mathcal{P}}

\def\cW{\mathcal{W}}

\def\fin{{\rm fin}}
\def\adm{{\rm sh}}
\def\re{{\rm re}}
\def\im{{\rm im}}
\def\nek{\text{\hbox{$\simeq$ \kern-.95em \hbox{$/$ \kern.05em}}}}

\def\opp{\operatornamewithlimits{\oplus}}

\DeclareMathOperator{\ad}{ad}

\DeclareMathOperator{\Id}{Id}

\DeclareMathOperator{\Supp}{Supp}

\DeclareMathOperator{\op}{opp}

\def\cplus{\hbox{$\subset${\raise1.05pt\hbox{\kern -0.55em
${\scriptscriptstyle +}$}}\ }}
\def\bcplus{\hbox{$\supset${\raise1.05pt\hbox{\kern -0.55em
${\scriptscriptstyle +}$}}\ }}

\def\ctimes{\hbox{$\times${\raise1.1pt\hbox{\kern -0.27em
${\scriptscriptstyle |}$}}\ }}

\def\udarrow{\hbox{$\nearrow${\kern -0.97em$\searrow$}\ }}

\def\bctimes{\hbox{$\times${\raise1.1pt\hbox{\kern -.74em
${\scriptscriptstyle |}$}}\ }\,\,}

\newtheorem{theorem}[equation]{Theorem}

\newtheorem{lemma}[equation]{Lemma}

\newtheorem{corollary}[equation]{Corollary}

\newtheorem{proposition}[equation]{Proposition}

\newtheorem{remark}[equation]{Remark}

\title{Classification of simple weight modules over affine Lie algebras}

\author{Ivan Dimitrov}
\thanks{Research supported by NSERC} 
\address{Department of Mathematics and Statistics, Queen's University, Kingston, K7L 3N6, Canada}
\email{dimitrov@mast.queensu.ca}
\author{Dimitar Grantcharov}
\address{Department of Mathematics, University of Texas at Arlington, Arlington, TX 76019, USA}
\email{grandim@uta.edu}

\begin{document}
\subjclass[2000]{Primary 17B10, Secondary 17B67}
\maketitle

\begin{abstract}
All simple weight modules with finite dimensional weight spaces over affine Lie algebras are classified.
\end{abstract}

\section*{Introduction} \label{secIntroduction}
The purpose of this paper is to complete the classification of all simple weight modules with finite dimensional
weight spaces over affine Lie algebras. This is a natural class of representations which contains all integrable
highest weight modules as well as the integrable loop modules with finite dimensional weight spaces. 
The integrable highest weight modules were the first class of representations over affine Lie algebras to be extensively
studied, see \cite{K} for detailed discussion of results and further bibliography. In \cite{C}  Chari 
classified all simple integrable weight modules with finite dimensional weight spaces over the untwisted affine Lie algebras. 
V. Chari and A. Pressley, \cite{CP2}, then extended this classification to all affine Lie algebras. The results of \cite{C} and \cite{CP2}
state that every simple integrable weight module is either a highest weight module or a loop module; in particular, the study of loop 
modules was initiated. In the 1990's V. Futorny began a comprehensive program studying weight modules over arbitrary affine Lie 
algebras. Here are just a few of the major results of this program: in \cite{F1} all parabolic sets of roots are described; the corresponding
parabolically induced modules are studied in \cite{F2}, \cite{F3}, and \cite{F6}; in \cite{F4} all simple weight modules of nonzero level 
with finite dimensional weight spaces are classified; in \cite{FKM} categories of parabolically induced modules are studied.
For a nice exposition of Futorny's program, see \cite{F5}. Still, the problem of complete classification was beyond reach until
O. Mathieu, \cite{M}, classified all simple weight modules with finite dimensional weight spaces over simple finite dimensional Lie algebras.
Following Mathieu's work Futorny in \cite{F5} and Futorny and Eswara Rao in \cite{EF} stated a number of conjectures whose
solution would be a major step towards solving the general classification problem. Finally, we mention that S. Berman and 
Y. Billig, \cite{BB}, constructed a class of simple parabolically induced modules with finite dimensional weight spaces, 
the so called exp--polynomial modules. 

We obtain our main result in three steps. First we establish a parabolic induction theorem which effectively proves Futorny's conjecture
that every cuspidal weight module, i.e. every simple weight module which is not parabolically induced,
 is dense for modules with finite dimensional weight spaces. Next we use an analog of Mathieu's
localization functor to establish that every cuspidal weight module is the twisted localization of a simple parabolically induced module.
Finally, we study the parabolically induced modules to describe what twisted localizations can be applied and to identify the resulting cuspidal 
modules. This scheme is very similar to the one followed by Mathieu. The major difference is that in step three we make extensive use of 
the constructions of Berman, Billig, Chari, and Pressley, as well as of Mathieu's classification to solve the difficult problem of identifying the 
simple constituents of twisted localizations of parabolically induced modules.

The paper is organized as follows. Section 1 contains background material. In Section 2 we establish a parabolic induction theorem which
reduces the classification problem to classifying all cuspidal weight modules. We also 
prove that every weight module with bounded weight multiplicities has finite length. 
In Section 3 we prove that every cuspidal weight 
module is the twisted localization of a simple parabolically induced module. We also provide an example of a simple graded module over
an affine Lie algebra which remains simple when considered without the grading.
In Section 4 we prove the main theorem which provides the complete list of cuspidal weight modules over affine Lie algebras. In particular,
we show that the twisted affine Lie algebras do not admit any cuspidal weight modules with finite dimensional weight spaces.


\medskip
\noindent
{\bf Acknowledgement.} We thank V. Futorny and J. Greenstein for useful discussions.

\medskip

\noindent
{\bf Notation.} 
The ground field $\F$ is algebraically closed of characteristic zero.
We denote the nonnegative real number by $\R_+$ and the nonnegative integers by $\Z_+$. 
The linear span of a subset $X$ of a vector space over a monoid $S$ is denoted by $\langle X \rangle_S$.

\section{Background results} \label{secNotation}

\medskip

\subsection{Categories of weight modules} \label{subsec_Weight}

Let $\gG$ be a Lie algebra and $\gH$ be a self--normalizing abelian subalgebra acting semisimply on $\gG$. Then $\gG$
decomposes as
\begin{equation} \label{eq0.1}
\gG = \gH \oplus (\opp_{\alpha \in \gH^*} \gG^\alpha),
\end{equation}
where $\gG^\alpha = \{x \in \gG \, | \, [h, x] = \alpha(h) x {\text { for every }} h \in \gH^*\}$. The subalgebra $\gH$ is called
a Cartan subalgebra 
of $\gG$, (\ref{eq0.1}) --- root decomposition of $\gG$ with respect to $\gH$, and 
$$
\Delta = \{ \alpha \in \gH^*\backslash\{0\} \, | \, \gG^\alpha \neq 0\}
$$
is the set of root of $\gG$ with respect to $\gH$. Denote by $\Delta^\re$ the set of roots $\alpha \in \Delta$ for which
$\Q \alpha \cap \Delta = \{ \pm \alpha \}$ and the subalgebra of $\gG$ generated by $\gG^\alpha + \gG^{-\alpha}$ is
isomorphic to $\gs\gl_2$. If $\alpha \in \Delta^\re$, then both $\gG^\alpha$ and  $\gG^{-\alpha}$ are one dimensional and
we fix nonzero elements $e_\alpha \in \gG^\alpha$ and $f_\alpha \in \gG^{-\alpha}$, which we refer to as root elements
corresponding to $\alpha$ and $-\alpha$ respectively. We also assume that $e_{-\alpha} = f_\alpha$. \footnote {This terminology 
is not standard but is consistent with the case when $\gG$ is a Kac--Moody Lie algebra.}

A $\gG$--module $M$ is a {\it weight module} if it is
$\gH$--diagonalizable, i.e. if
$$
M = \oplus_{\lambda \in \gH^*} M^\lambda, {\text { where }} M^\lambda = \{ m \in M \, | \, h \cdot m = \lambda(h) m {\text { for every }} h \in \gH\}.
$$
The {\it support of} $M$ is the set $\Supp M := \{ \lambda \in \gH^* \, | \, M^\lambda \neq 0\}$. Denote the category of weight $\gG$--modules
by $\cW$ (when the Lie algebra whose modules we consider is not clear, we will write $\cW(\gG)$ instead).
The full subcategory of $\cW$ consisting of $\gG$--modules with finite dimensional weight spaces
is denoted by $\cW_\fin$. 
If $\alpha \in \Delta$, a module $M \in \cW_\fin$ is {\it $\alpha$--bounded} if $\sup_{q \in \Q} \dim M^{\lambda + q \alpha} < \infty$ for
every $\lambda \in \Supp M$.

A module $M \in \cW_\fin$  {\it has a shadow} if for every $\alpha \in \Delta^{\re}$ the root element $e_\alpha$ acts either locally nilpotently or injectively on $M$.
The full subcategory of $\cW_\fin$ consisting of all $\gG$--modules with shadow is denoted by $\cW_\adm$.
 and the full subcategory of
$\cW_\fin^\Phi$ consisting of modules with a shadow is denoted by $\cW_\adm^\Phi$.

Every $M \in \cW_\adm$ determines
a decomposition of $\Delta^\re$ as $\Delta^\re = \Delta^f(M) \sqcup \Delta^i(M)$, where $\Delta^f(M)$ and $\Delta^i(M)$ consist of
 the roots whose corresponding
root elements act locally nilpotently and injectively on $M$ respectively. If the module $M$ is clear from the context we may write 
$\Delta^f$ and $\Delta^i$ instead of $\Delta^f(M)$ and $\Delta^i(M)$. The following properties of modules from $\cW_\adm$ are
straightforward. 

\medskip

\begin{proposition} \label{prop3} If $M \in \cW_\fin$ is simple then $M \in \cW_\adm$. If $M \in \cW_\adm$ then the following hold.
\begin{itemize}
\item[\rm{(i)}] If $N$ is a nontrivial submodule of $M$, then $N \in \cW_\adm$ and
$\Delta^f(N) = \Delta^f(M)$ and $\Delta^i(N) = \Delta^i(M)$.
\item[\rm{(ii)}] If $\alpha \in \Delta^\re$ then the following are equivalent
\begin{itemize}
\item[\rm{(a)}] $\alpha \in \Delta^f$;
\item [\rm{(b)}]  $\{\lambda_0 + n \alpha \, | \, n \in \Z_+\} \cap \Supp M$
is a finite set for some $\lambda_0 \in \Supp M$;
\item [\rm{(c)}]   $\{\lambda + n \alpha \, | \, n \in \Z_+\} \cap \Supp M$  is a finite
set for every $\lambda \in \Supp M$.
\end{itemize}
\item[\rm{(iii)}] If $\alpha \in \Delta^\re$ then the following are equivalent
\begin{itemize}
\item[\rm{(a)}] $\alpha \in \Delta^i$;
\item [\rm{(b)}]  $\{\lambda_0 + n \alpha \, | \, n \in \Z_+\} \cap \Supp M$
is an  infinite set for some $\lambda_0 \in \Supp M$;
\item [\rm{(c)}]   $\{\lambda + n \alpha \, | \, n \in \Z_+\} \cap \Supp M$  is an infinite
set for every $\lambda \in \Supp M$.
\end{itemize}
\item[\rm{(iv)}] If  $-\alpha \in \Delta^i$, then the following are equivalent
\begin{itemize}
\item[\rm{(a)}] $\alpha \in \Delta^i$;
\item[\rm{(b)}] $f_\alpha$ is bijective;
\item[\rm{(c)}] $\dim M^\lambda = \dim M^{\lambda + n \alpha}$ for every $\lambda \in \Supp M$ and every $n \in \Z$.
\end{itemize}
\end{itemize}
\end{proposition}

\medskip

\subsection{Affine Lie algebras} \label{subsecAffine}
We recall the the construction of affine Lie algebras and fix notation. For more detail, see \cite{K}.

Let $\gg$ be a simple finite dimensional Lie algebra with  non--degenerate invariant symmetric bilinear form 
$(\, \, , \, ): \gg \times \gg \to \F$. (Usually, $(\, \, , \, )$ will be the Killing form of $\gg$.) 
Denote by $\cL(\gg)$ the loop algebra $\gg \otimes \F[t, t^{-1}]$. The affine Lie algebra 
$\cA(\gg) = \cL(\gg) \oplus \F D \oplus \F K$ has commutation relations
$$
[x \otimes t^m, y \otimes t^n] = [x,y] \otimes t^{m+n} + \delta_{m, -n} m \, (x,y) K, \quad
[D, x \otimes t^m] = m x \otimes t^m,  \quad [K, \cA(\gg)] =0,
$$
where $x, y \in \gg$, $m, n \in \Z$, and $\delta_{i,j}$ is Kronecker's delta. The form $(\, \, , \, )$ extends to 
a non--degenerate invariant symmetric bilinear form on $\cA(\gg)$, still denoted by 
$(\, \, , \, ): \cA(\gg) \times \cA(\gg) \to \F$, via
\begin{align*}
& (x \otimes t^m, y \otimes t^n) = \delta_{m, -n} (x, y), \quad (D, K) = 1, \\
& (x \otimes t^m, D) = (x \otimes t^m, K) = (K, K) = (D, D) = 0.
\end{align*}

If $\sigma$ is a diagram automorphism of $\gg$ of order $s$ then $\sigma$ extends to an automorphism of $\cA(\gg)$ by 
$$
\sigma( x \otimes t^m) = \zeta^m \sigma(x) t^m, \quad \sigma(D) = D, \quad \sigma(K) = K,
$$
 where $\zeta$ is a fixed primitive 
 $s^{th}$ root of unity. The twisted affine Lie algebra $\cA(\gg, \sigma)$ is  the Lie algebra $\cA(\gg)^\sigma$ of
 $\sigma$--fixed points of $\cA(\gg)$. Note that 
 $$
 \cA(\gg,\sigma) = \cL(\gg,\sigma) \oplus \F D \oplus \F K,
 $$
 where $\cL(\gg, \sigma) = \cL(\gg)^\sigma$ is the subalgebra of $\sigma$--fixed points of $\cL(\gg)$. 
 One checks immediately that 
 $$
\cL(\gg, \sigma) = \opp_{j \in \Z} \gg_{\bar{j}} \otimes t^j,
$$
 where 
 $$
 \gg = \opp_{\bar{j} \in \Z/s\Z} \gg_{\bar j}
 $$
 is the decomposition of $\gg$ into $\sigma$--eigenspaces.
 The restriction
 of $(\, \, , \, )$ to $\cA(\gg, \sigma)$ is a non--degenerate invariant symmetric bilinear form for which we will use the same notation.
 
If $\gh$ is a Cartan subalgebra of $\gg$, then $\gh \oplus \F D \oplus \F K$ is a Cartan subalgebra of $\cA(\gg)$.
Furthermore, $\sigma$ preserves $\gh$ and $\gh^\sigma \oplus \F D \oplus \F K$ is a Cartan subalgebra of $\cA(\gg, \sigma)$.
For the rest of the paper $\gh$ denotes a fixed Cartan subalgebra of $\gg$, $\gG$ denotes $\cA(\gg)$ or $\cA(\gg,\sigma)$,
 and $\gH$ denotes  the corresponding Cartan subalgebra of $\gG$. \footnote{One can consider the cases of $\cA(\gg)$ and 
 $\cA(\gg, \sigma)$ simultaneously by taking $\sigma$ to be the identity in the case of $\cA(\gg)$. We will distinguish between the two cases in the hope that this improves the clarity of the text.}

The Lie algebra $\gG$ admits a root decomposition
$$
\gG = \gH \oplus (\opp_{\alpha \in \Delta} \gG^\alpha).
$$
To describe the root system $\Delta$ of $\gG$, let $\delta \in \gH^*$ denote the element with
$$
\delta(D) = 1, \quad \delta(x \otimes t^m) = \delta(K) = 0.
$$
If $\gG = \cA(\gg)$, denote the root system of $\gg$ by $\mathring{\Delta}$. If $\gG = \cA(\gg, \sigma)$,
denote the nonzero weights of the $\gg_{\bar{0}}$--module $\gg_{\bar{j}}$ by $\mathring{\Delta}_{\bar{j}}$ and set
$\mathring{\Delta} = \cup_{\bar{j} \in \Z/s\Z} \mathring{\Delta}_{\bar{j}}$. Note that for 
$\gG = \cA(\gg, \sigma) \not \cong A_{2l}^{(2)}$, $\mathring{\Delta} = \mathring{\Delta}_{\bar{0}}$ is the root system of
$\gg_{\bar{0}}$ and for $\gG \cong A_{2l}^{(2)}$, $\mathring{\Delta}$ is the non--reduced root system $BC_l$.

The decomposition $\gH = \gh \oplus \F D \oplus \F K$
(respectively, $\gH = \gh^\sigma \oplus \F D \oplus \F K$) allows us to consider $\mathring{\Delta}$ as a subset of $\gH^*$.
The root system $\Delta$ decomposes as 
$$
\Delta = \Delta^{\re} \sqcup \Delta^{\im},
$$
 where 
 $$
 \Delta^{\im} = \{n \delta \, | \, n \in \Z \backslash \{0\}\}
 $$ 
 are the imaginary roots of $\gG$ and the real roots $\Delta^{\re}$ are given as follows.
\begin{itemize}
\item[(i)] If $\gG = \cA(\gg)$, then
$$
\Delta^\re = \{ \alpha + n \delta \, | \, n \in \Z, \, \alpha \in \mathring{\Delta}\}.
$$
\item[(ii)] If $\gG = \cA(\gg, \sigma)$, then
$$
\Delta^\re = \{ \alpha + n \delta \, | \, n \in \Z, \, \alpha \in \mathring{\Delta}_{\bar{n}}\}.
$$
\end{itemize}

For every real root $\alpha\in \Delta^\re$, $r \alpha \in \Delta$ with $r \in \R$ if and only if $r = \pm 1$. Moreover, $\dim \gG^\alpha = 1$ and there exist elements $e_\alpha \in \gG^\alpha$,
$f_\alpha \in \gG^{-\alpha}$, and $h_\alpha \in \gH$, such that $(e_\alpha, h_\alpha, f_\alpha)$ is an $\gs\gl(2)$--triple. For
the rest of the paper we fix such elements $e_\alpha, h_\alpha, f_\alpha$  for every $\alpha \in \Delta^{\re}$ and assume that $e_{-\alpha} = - f_\alpha$, $h_{-\alpha} = - h_\alpha$,
and $f_{-\alpha} = - e_{\alpha}$.

A subset $\Phi$ of $\Delta$ is called a base of $\gG$ if every element of $\Delta$ can be written uniquely as an integer
combination of elements of $\Phi$ in which either all coefficients are non--negative or all coefficients are non--positive. 
The elements of $\Phi$ are called simple roots, 
the roots which are expressed as non--negative combination of simple roots --- positive roots, and the roots opposite 
to positive roots --- negative roots. 
Every affine Lie algebra admits a base, all bases are conjugate under the action of $W \times \Z_2$, and every base
consists of real roots. 
 
We will also need the definition of Heisenberg algebra. It is the Lie algebra $\cH$ with basis $\{X_n, D, K \, | \, n \in \Z\}$ and commutation relations
$$
[X_m, X_n] = \delta_{m,-n} m K, \quad \quad [D, X_n] = n X_n, \quad \quad [K, \cH] = 0.
$$
The subalgebra $\F X_0 \oplus \F D \oplus \F K$ is a Cartan subalgebra with respect to which $\cH$ has a root decomposition with root system
$\{n \delta \, | \, n \in \Z \backslash \{0\}\}$. For uniformity of notation, we write $\cH = \cL(\F X_0) \oplus \F D \oplus \F K$, i.e. we consider $\cH$ 
as the affine
Lie algebra corresponding to a one dimensional (abelian) Lie algebra.

\medskip

\subsection{Parabolic subalgebras} \label{subsecParabolic}
Following Bourbaki, \cite{Bo}, we call a subset $P \subset \Delta$ parabolic if 
\begin{itemize}
\item[(i)] $\Delta = P \cup -P$,
\item[(ii)] $\alpha, \beta \in P$ with $\alpha + \beta \in \Delta$ implies $\alpha + \beta \in P$.
\end{itemize}
A parabolic subset $P$ defines a parabolic
subalgebra $\gP$ of $\gG$ by setting 
\begin{equation} \label{eq1.1}
\gP = \gH \oplus (\opp_{\alpha \in P} \gG^\alpha).
\end{equation}
The subalgebras 
\begin{equation} \label{eq1.11}
\gL = \gH \oplus (\oplus_{\alpha \in P \cap (-P)} \gG^\alpha) \quad 
{\text {and }} \quad \gN^+ = \oplus_{\alpha \in P \backslash (-P)} \gG^\alpha
\end{equation}
are the Levi component and the nillradical of $\gP$ respectivelly.

A parabolic subalgebra for which $P \cap -P = \emptyset$ is called a Borel subalgebra. If $\Phi$ is a base of $\Delta$, 
the set of positive roots is a parabolic subset of $\Delta$.

The parabolic subsets of roots of affine Lie algebras were described explicitly by V. Futorny in \cite{F1}.
Below we provide this classification in slightly different terms. For details, see \cite{DFG}.

Denote by $Q$ the root lattice of $\gG$, i.e. the abelian group generated by $\Delta$. Let $V = \R \otimes_\Z Q$ be the real vector space spanned by $\Delta$. 
Every linear function $\varphi : V \to \R$ determines a decomposition
\begin{equation} \label{eq1.2}
\Delta = \Delta^+ \sqcup \Delta^0 \sqcup \Delta^-,
\end{equation}
where $\Delta^\pm = \{\alpha \in \Delta \, | \, \varphi(\alpha) \gtrless 0\}$ and $\Delta^0 = \{\alpha \in \Delta \, | \, \varphi(\alpha) = 0\}$. We call a decomposition of the form (\ref{eq1.2}) 
a triangular decomposition of $\Delta$. Clearly, the set $\Delta^+ \sqcup \Delta^0$ is a parabolic subset of $\Delta$. A parabolic subset of $\Delta$ of the form 
(\ref{eq1.2}) is called principal parabolic subset. It is not true, however, that every parabolic subset of $\Delta$ is principal. Indeed, if 
$$
\Delta^0 = (\Delta^0)^+ \sqcup (\Delta^0)^0 \sqcup (\Delta^0)^-
$$
is a triangular decomposition of $\Delta^0$, then
\begin{equation} \label{eq1.3}
P = \Delta^+ \sqcup (\Delta^0)^0 \sqcup (\Delta^0)^+
\end{equation}
is also a parabolic subset of $\Delta$ which is not necessarily principal. Theorem 4.13 of \cite{DFG} implies 
that every parabolic subset of $\Delta$ is of the form (\ref{eq1.3}). For convenience we provide the statement of
Theorem 4.13 of \cite{DFG} below.

\medskip

\begin{proposition} \label{prop_parabolic} 
If  $P \subset \Delta$ is a parabolic set of roots, then $P$ is of the form (\ref{eq1.3}).
Furthermore, if $P$ is a proper subset of $\Delta$, one of the three mutually exclusive alternatives holds.

\begin{itemize}
\item[(i)] $P$ is principal and $\delta \not \in \Delta^0$;

\item[(ii)] $P$ is principal and $(\Z \setminus \{0\}) \delta \subset P$;

\item[(iii)] $P$ is not principal.
\end{itemize}
\end{proposition}

\medskip

The subalgebras (\ref{eq1.1}) defined by parabolic sets from cases (i), (ii), and (iii) of Proposition \ref{prop_parabolic} are called standard, imaginary, and mixed type parabolic subalgebras of $\gG$ respectively. We also use the terms standard, imaginary,
and mixed type Borel subalgebras when $P \cap (-P) = \emptyset$.

It is not difficult to notice that the triangular decomposition (\ref{eq1.2}) is not uniquely determined by $P$. However, there is a unique triangular decomposition
(\ref{eq1.2}) corresponding to $P$ for which $\Delta^0$ is minimal. For the rest of the paper by a triangular decomposition of $\Delta$ corresponding to $P$ we will
mean the one with minimal $\Delta^0$. It turns out that standard and imaginary parabolic subalgebras of $\cG$ give rise to simple weight $\gG$--modules 
with finite dimensional weight spaces, while mixed type parabolic subalgebras do not. For this reason we will concentrate on studying standard and imaginary 
parabolic subalgebras corresponding to principal triangular decompositions (\ref{eq1.2}).

If $\gB$ is a Borel subalgebra of $\gG$ then $\gB$ is standard if and only if $\gB$ admits a base, i.e., there exists a base
$\Phi$ of $\Delta$ for which $P$ is the corresponding set of positive roots. Furthermore, a parabolic subalgebra $\gP$ is
standard if and only if it contains a standard Borel subalgebra.

Consider a triangular decomposition (\ref{eq1.2}) of $\Delta$, the corresponding parabolic subalgebra $\gP$,
the Levi component $\gL$ of $\gP$, and the nillradical $\gN^+$ of $\gP$, see (\ref{eq1.1}) and (\ref{eq1.11}).

If $\gP$ is standard, then $\Delta^0$ is a finite root system and $\gL$ is a finite dimensional reductive Lie algebra.
If $\gP$ is imaginary, then
 $\Delta^0$ is the union of root systems of affine Lie algebras, possibly equal to $\Delta^\im$. The algebra $\gL$ in this case is a
subfactor of a direct sum $\gL_1 \oplus \gL_2 \oplus \ldots \oplus \gL_k$, where each $\gL_i$ is isomorphic to an affine Lie algebra or to the Heisenberg algebra.
More precisely, let $\gL_i = \cL(\gg_i, \sigma_i) \oplus \F D_i \oplus \F K_i$. Then there exists $D \in \F D_1 \oplus \F D_2 \oplus \ldots \oplus \F D_k$, such that
$$
\gL \simeq (\opp_{i = 1}^k \cL(\gg_k, \sigma_k) ) \oplus \F D \oplus \F K,
$$
where $K = K_1=K_2= \ldots = K_k$.

\medskip

\subsection{Parabolically induced modules} \label{subsec_1.5}
Let $\gP$ be a parabolic subalgebra of $\gG$ with Levi component $\gL$ and nillradical $\gN^+$.
Recall that $Q$ denotes the root lattice of $\gG$ and denote the root lattice of $\gL$ by $Q^0$.
If $N$ is a weight $\gL$--module whose support 
is contained in a single $Q^0$--coset $\lambda + Q^0 \in
\gH^*/Q^0$, then $N$ can be endowed with a $\gP$--module structure  
 with the trivial action of $\gN^{+}$. Set
$$
M_{\gP}(N):=U(\gG) \otimes_{U(\gP)} N.
$$
Among the submodules of $M_\gP(N)$ which intersect trivially with $N$ 
there is a unique maximal element $Z_{\gP}(N)$. Set 
$$
V_\gP(N) = M_{\gP}(N)/Z_{\gP}(N).
$$
The following proposition is straghtforward.

\medskip

\begin{proposition} \label{ind_simple}
Let $N$ be a weight $\gL$--module whose support is contained in a
single $Q^0$--coset. Then

\noindent 
\rm{(i)} $(V_{\gP}(N))^{\gN^+}= N$;

\noindent 
\rm{(ii)}  $V_{\gP}(N)$ is a simple $\gG$--module if and only if $N$ is simple $\gL$--module;

\noindent
\rm{(iii)} If $M$ is a simple $\gG$--module with $M^{\gN^+} \neq 0$ then $M \cong V_\gP(M^{\gN^+})$;

\noindent
\rm{(iv)} If $\gP$ is a standard parabolic subalgebra of $\gG$ then 
$V_{\gP}(N) \in \cW_\fin(\gG)$ if and only if $N \in \cW_\fin(\gL)$.
\end{proposition}

\medskip

Next we recall a result due to Berman and Billig, \cite{BB}.

\medskip

\begin{proposition} \label{prop_exp-pol}
Let $\gP$ e a parabolic subalgebra of $\gG$ with Levi component $\gL$. Assume that $N$ is a simple weight $\gL$--module 
with finite dimensional weight spaces which 
admits a basis in which the coefficients of the action of a root basis of $\gL$ are exp--polynomial functions.
Then $V_\cP(N)$ has  finite dimensional weight spaces.
\end{proposition}

\medskip

\subsection{Loop modules} \label{subsee1.55} 
We recall the definition and some properties of loop modules. For more detail, see \cite{C}, \cite{CP1}, \cite{CP2}.

Let $\gG = \cA(\gg)$, let $V_1, \ldots, V_k$ be weight $\gg$--modules, and let $a_1, \ldots, a_k$ be nonzero scalars. 
Following Chari and Pressley we define the loop module $\cL_{a_1, \ldots, a_k}(V_1 \otimes \ldots \otimes V_k)$ in the following
way: the underlining vector space of $\cL_{a_1, \ldots, a_k}(V_1 \otimes \ldots \otimes V_k)$ is 
$$(V_1 \otimes \ldots \otimes V_k) \otimes \F[t,t^{-1}],$$
$X \otimes t^n \in \cA(\gg)$ acts as
$$
(X \otimes t^n) \cdot ((v_1 \otimes \ldots \otimes v_k) \otimes t^s) = \sum_{i=1}^k a_i^n (v_1 \otimes \ldots \otimes X\cdot v_i \otimes \ldots \otimes v_k)
\otimes t^{n+s},
$$
$D$ acts as $t\frac{d}{dt}$, and $K$ acts trivially. If the scalars $a_1, \ldots, a_k$ are distinct, then $\cL_{a_1, \ldots, a_k}(V_1 \otimes \ldots \otimes V_k)$ 
is completely reducible and decomposes as a direct sum of finitely many isomorphic simple $\cL(\gg)$--modules. Denote by $V_{a_1, \ldots, a_k}(V_1 \otimes \ldots \otimes V_k)$ the simple $\cL(\gg)$--module which is a component of $\cL_{a_1, \ldots, a_k}(V_1 \otimes \ldots \otimes V_k)$.
Considered as $\cA(\gg)$--modules, 
the constituents of $\cL_{a_1, \ldots, a_k}(V_1 \otimes \ldots \otimes V_k)$ differ only by a shift of the action of $D$.
By a slight abuse of notation we denote by $V_{a_1, \ldots, a_k}(V_1 \otimes \ldots \otimes V_k)$ any shift of a simple $\cA(\gg)$--component
of $\cL_{a_1, \ldots, a_k}(V_1 \otimes \ldots \otimes V_k)$.

If $\gG = \cA(\gg, \sigma)$, then $\cL_{a_1, \ldots, a_k}(V_1 \otimes \ldots \otimes V_k)$ admits an endomorphism compatible with $\sigma$
if and only if the modules $V_1, \ldots, V_k$ come in $r$--tuples of isomorphic modules and for each $r$--tuple the corresponding scalars
$a_1, \ldots, a_r$ are multiples (with the same scalar) of the $r^{th}$ roots of unity. We denote the corresponding endomorphism of
$\cL_{a_1, \ldots, a_k}(V_1 \otimes \ldots \otimes V_k)$ by $\sigma$ again and let $\cL_{a_1, \ldots, a_k}^\sigma(V_1 \otimes \ldots \otimes V_k)$
denote the fixed points of $\sigma$. Similarly, we can define $V_{a_1, \ldots, a_k}^\sigma(V_1 \otimes \ldots \otimes V_k)$.

\subsection{Localization. Elementary properties} \label{subsec1.6}
Let us first recall the definition of the localization functor of weight modules.
For details we refer the reader to \cite{De} and \cite{M}.

Denote by $U$ the universal enveloping algebra $U(\gG)$ of $\gG$.
For every $\alpha \in \Delta^\re$ the multiplicative set
$F_{\alpha}:=\{ f_\alpha^n \; | \; n \in \Z_+ \} \subset U$ 
satisfies Ore's localization conditions
because $\ad f_\alpha$ acts locally nilpotently on $U$. Let
$U_\alpha$ be the localization of $U$ relative to $F_{\alpha}$.
For every $M \in \cW$ we denote by $\cD_\alpha M$ the {\it
$\alpha$--localization of $M$}, defined as $\cD_\alpha M =
U_\alpha \otimes_U M$. If $f_\alpha$ is injective on $M$, then $M$
is a submodule of $\cD_\alpha M$, $f_\alpha$ is injective on
$D_\alpha M$, and $D_\alpha^2 M = D_\alpha M$. Furthermore, if
$f_\alpha$ is injective on $M$, then it is bijective on $M$ if and
only if $\cD_\alpha M = M$, cf. Proposition \ref{prop3}. Finally,
if $[f_\alpha, f_\beta] = 0$ and both $f_\alpha$ and $f_\beta$ are
injective on  $M$, then $\cD_\alpha \cD_\beta M = \cD_\beta
\cD_\alpha M$.

In what follows we recall the definition of a generalized conjugation in $U_\alpha$ introduced in \cite{M}.
 For $x \in \F$ and $u \in U_\alpha$  we set
\begin{equation} \label{theta}
\Theta_x(u):= \sum\limits_{i=0}^\infty \binom{x}{i}\,
\ad(f_\alpha)^i (u) \, f_\alpha^{-i},
\end{equation}
where $\binom{x}{i}= \frac{x(x-1)...(x-i+1)}{i!}$. Since $\ad
(f_\alpha)$ is locally nilpotent on $U_\alpha$, the sum above is
actually finite. Note that for $x \in \Z$ we have $\Theta_x(u) =
f_\alpha^x u f_\alpha^{-x}$.  For a $U_\alpha$-module $M$ by
$\Phi_\alpha^x M$ we denote the $U_\alpha$-module $M$ twisted by
the action
$$
u \cdot v^x := ( \Theta_x (u)\cdot v)^x,
$$
where $u \in U_\alpha$, $v \in M$, and $v^x$ stands for the
element $v$ considered as an element of $\Phi_\alpha^x M$. In
particular, $v^x \in M^{\lambda + x \alpha}$ whenever $v \in
M^\lambda$. Since $v^n = f_{\alpha}^{-n} \cdot v$ whenever $n \in
\Z$ it is convenient to set $f_{\alpha}^x \cdot v :=v^{-x}$ for $x
\in \F$.

The following lemma is straightforward.

\medskip

\begin{lemma} \label{lmnew} In the category of all $U_{\alpha}$-modules  we have

\noindent {\rm (i)} $\Phi_\alpha^x  = \Id$ whenever $x \in \Z$.

\noindent {\rm(ii)} $\Phi_\alpha^x \circ
\Phi_\alpha^y=\Phi_\alpha^{x+y}$ and, consequently, $\Phi_\alpha^x
\circ\Phi_\alpha^{-x} = \Phi_\alpha^{-x} \circ\Phi_\alpha^{x}  =
\Id$ for any $x,y \in \F$.

\noindent {\rm(iii)} $f_{\alpha}^x \cdot (f_{\alpha}^y \cdot v) =
f_{\alpha}^{x+y} \cdot v$ for any $x,y \in \F$.

\noindent {\rm(iv)} $f_{\alpha}^x \cdot (u \cdot (f_{\alpha}^{-x}
\cdot v)) = \Theta_x(u) \cdot v$ for any $x \in \F$.

\end{lemma}

\medskip

\begin{remark} \label{rem_loc}
We will often consider $\cD_{\alpha}$ and $\Phi_{\alpha}^x$ over
more general categories of modules. Note that, $\cD_{\alpha} M$ is
well-defined for any 
$\gG^{\alpha}$--module $M$, and $\Phi_{\alpha}^x N$ is
well-defined for any $U(\gG^{\alpha})_{\alpha}$--module $N$. Furthermore, if $M$ is $f_{\alpha}$--injective,
then $M \to \cD_{\alpha} M$ is a monomorphism.
\end{remark}

\medskip

In what follows we set $\cD_\alpha^x M:=\Phi_\alpha^x (\cD_\alpha M)$ and refer to it as a twisted localization of $M$.

\medskip

\section{Parabolic induction theorem} \label{sec_Induction}

\medskip

\subsection{The roots of the nillradical of a standard $\gP$} \label{subsec_cone}
If  $\gP$ is a  standard parabolic subalgebra of $\gG$ corresponding to the triangular decomposition (\ref{eq1.2}),
then 
$\Delta^0$ is a finite root system and $\gL$ is a finite dimensional reductive
Lie algebra. Fix a Borel subalgebra $\gB$ of $\gG$ contained in $\gP$ with base $\Phi = \{\alpha_0, \ldots, \alpha_l\}$.
Assume that $\alpha_i \in \Delta^0$ if and only if $i \in I \subset \{0, \ldots, l\}$ and set $J = \{0, \ldots, l\} \backslash I$.
Let
\begin{equation} \label{eq1.23}
\Phi_\gP = \{w(\alpha_j)\}_{w \in W_\gL, j \in J},
\end{equation}
where $W_\gL$ denotes the Weyl group of $\gL$. Clearly $\Phi_\gP$ does not depend on the choice of 
$\gB \subset \gP$.
Let, furthermore,  $\cC_\gP^+ = \langle \Phi_\gP \rangle_{\Q_+}$ denote the cone generated by $\Phi_\gP$, 
 considered  as a subset of the
$\Q$--vector space $\langle \Delta \rangle_{\Q}$. Finally, denote by $Q_\gP$ the abelian group generated by $\Phi_\gP$.

\medskip

\begin{proposition} \label{prop_cone} For any standard parabolic subalgebra $\gP \subset \gG$ with
$\delta \in \Delta^+$ the imaginary root $\delta$ belongs to the interior of $\cC_\gP^+$. 
The group $Q_\gP$ has a finite index in $Q$. Moreover there exists a constant $N_\gG$
depending on $\gG$ only such that for every $\nu \in Q$ and every standard parabolic subalgebra $\gP$ we have
$N_\gG \nu \in Q_\gP$.
\end{proposition}

\medskip

\noindent
{\bf Proof.}  If $\delta \in \Delta^+$, then 
$$
\delta = c_0 \alpha_0 + c_1 \alpha_1 + \ldots + c_l \alpha_l
$$
with $c_i > 0$ for $0 \leq i \leq l$. 
Since $\delta$ is $W$--invariant (and, in particular, $W_{\gL}$--invariant) we obtain
$$
|W_\gL| \delta = \sum_{w \in W_\gL} \sum_{i=0}^{l} c_i w(\alpha_i) = \sum_{w \in W_\gL} \sum_{i \in I} c_i w(\alpha_i) 
+ \sum_{w \in W_\gL} \sum_{j \in J} c_j w(\alpha_j),
$$
where $|W_\gL|$ as usual denotes the order of $W_\gL$. 
Note that $\sum_{w \in W_\gL} w(\alpha_i) = 0$ for every $i \in I$ because it is a $W_\gL$--invariant element of
$Q^0$. Thus
\begin{equation} \label{eq1.24}
|W_\gL| \delta = \sum_{\beta \in \Phi_\gP} d_\beta \beta \quad {\text { with }} \quad d_\beta > 0 \quad 
{\text { for every }} \quad \beta \in \Phi_\gP.
\end{equation}

Next we show that $\Phi_\gP$ spans $\langle \Delta \rangle_{\Q}$. Set
$$
I' = \{ i \in I \, | \, (\alpha_i, \beta) = 0 {\text { for every }} \beta \in \Phi_\gP\}
$$
and $I'' = I \backslash I'$. First we show that $I' = \emptyset$. Indeed, assume that $i \in I'$ and $j \in I'' \cup J$.
We claim that $(\alpha_i, \alpha_j) = 0$. If $j \in J$, then $\alpha_j \in \Phi_\gP$ and $(\alpha_i, \alpha_j) = 0$ from
the definition of $I'$. If, on the other hand, $j \in I''$, then there exists $\beta \in \Phi_\gP$ such that $(\alpha_j, \beta) 
\neq 0$. Note that $s_{\alpha_j}(\beta) \in \Phi_\gP$ and we have
$$
0 = (\alpha_i, s_{\alpha_j}(\beta)) = \left(\alpha_i, \beta - \frac{2 (\alpha_j, \beta)}{(\alpha_j, \alpha_j)} \alpha_j \right)=
(\alpha_i, \beta) - \frac{2 (\alpha_j, \beta)}{(\alpha_j, \alpha_j)}(\alpha_i, \alpha_j).
$$
Using the assumptions on $\alpha_i$ and $\alpha_j$ we conclude that $(\alpha_i, \alpha_j) = 0$. 
Thus we obtained that the sets $\{\alpha_i\}_{i \in I'}$ are $\{\alpha_j\}_{j \in I'' \cup J}$ are orthogonal. Since
$I'' \cup J \neq \emptyset$ we conclude that $I' = \emptyset$. In other words, for every $i \in I$ there exists
$\beta \in \Phi_\gP$ with $(\alpha_i, \beta) \neq 0$. Consequently,
\begin{equation} \label{eq1.25}
\alpha_i = \frac{(\alpha_i, \alpha_i)}{2 (\alpha_i, \beta)} (\beta - s_{\alpha_i}(\beta)) \in \langle \Phi_\gP \rangle_\Q,
\end{equation}
which concludes the proof that $\Phi_\gP$ spans $\langle \Delta \rangle_{\Q}$. Combining this with
 equation (\ref{eq1.24}) we obtain that $\delta$ belongs to the  interior of $\cC_\gP^+$. 

The abelian group $Q_\gP$ has a finite index in $Q$
because $\Phi_\gP$ and $\Delta$ generate the same $\Q$--vector space.
Moreover, equation (\ref{eq1.25}) shows that 
 if we denote the 
 least common multiple of the denominators of the fractions
$\frac{(\alpha_i, \alpha_i)}{2 (\alpha_i, \beta)}$ with $i \in I$ and $\beta \in \Phi_\gP$ such that  $(\alpha_i, \beta) \neq 0$
by $N_\gP$, then $N_\gP \nu \in Q_\gP$ for any $\nu \in Q$.

For $\beta = w(\alpha_j)$ with $w \in W_{\gL}$ and $j \in J$ we have
$$
\frac{(\alpha_i, \alpha_i)}{2 (\alpha_i, \beta)} = \frac{(\alpha_i, \alpha_i)}{2 (\alpha_i, w(\alpha_j))} =
\frac{(w^{-1}(\alpha_i), w^{-1}(\alpha_i))}{2 (w^{-1}(\alpha_i), \alpha_j)} = \frac{(\pi(w^{-1}(\alpha_i)), \pi( w^{-1}(\alpha_i)))}{2 (\pi( w^{-1}(\alpha_i)), \pi(\alpha_j))},
$$
where $\pi : \Delta \to \mathring{\Delta}$ is the natural projection. This equation
 shows that if we define 
$N_\gG$ as the least common multiple of the denominators of all fractions in the set
$$
\left\{\frac{(\alpha', \alpha')}{2 (\alpha', \beta`)} \, |\, \alpha`, \beta` \in  \mathring{\Delta} \right\},
$$
then $N_\gG \nu \in Q_\gP$ for every $\nu \in Q$ and every standard parabolic subalgebra $\gP \subset \gG$.
\hfill $\square$

\medskip

\subsection{Existence of extreme weights} \label{subsec_extreme}

\medskip

\begin{proposition} \label{prop2.0}
Let $\Delta = \Delta^- \sqcup \Delta^0 \sqcup \Delta^+$ be a triangular decomposition of $\Delta$ and let 
$M \in \cW_\adm$ be a $\gG$--module with $\Delta^+ \cap \Delta^\re \subset \Delta^f$ and $\Delta^- \cap \Delta^\re \subset \Delta^i$. There exists $\lambda \in \Supp M$ 
such that $\lambda + \beta \not \in \Supp M$ for every $\beta \in \Delta^+$.
\end{proposition}

\medskip

\noindent
{\bf Proof.} We consider separately the cases when $\delta \not \in \Delta^0$ and $\delta  \in \Delta^0$.

\noindent
{\bf Case 1: $\delta \not \in \Delta^0$.} We assume without loss of generality that $\delta \in \Delta^+$. The parabolic 
subalgebra $\gP$ is standard and Proposition \ref{prop_cone} applies. Using the notation of Proposition \ref{prop_cone}
we observe that $(\lambda + \cC_{\gP}^+) \cap \Supp M = \emptyset$ whenever $\lambda \not \in \Supp M$. Set
$\Phi_\gP = \{ \beta_1, \ldots, \beta_N\}$.

We are going to construct a sequence $\lambda_1, \ldots, \lambda_N$ of weights in $\Supp M$ such that
\begin{equation} \label{eq2.11}
(\lambda_i + \beta_j + \cC_{\gP}^+) \cap \Supp M = \emptyset 
\end{equation}
 for every $1\leq j \leq i \leq N$.
We proceed by induction.
Fix $\lambda \in \Supp M$. Since $\beta_1 \in \Delta^f$, there exists $\lambda_1 = \lambda + k_1 \beta_1 \in
\Supp M$ such that $\lambda_1 + \beta_1 \not \in \Supp M$. Hence $(\lambda_1 + \beta_1 + \cC_{\gP}^+) \cap \Supp M = \emptyset$.
Assume $\lambda_1, \ldots, \lambda_s$ have been chosen so that (\ref{eq2.11}) holds for every $1 \leq j \leq i \leq s$. 
We choose $\lambda_{s+1} = \lambda_s + k_{s+1} \beta_{s+1} \in \Supp M$ so that $\lambda_{s+1} + \beta_{s+1} \not \in \Supp M$. 
This implies that $(\lambda_{s+1} + \beta_{s+1} + \cC_{\gP}^+) \cap \Supp M = \emptyset$. Furthermore, 
$(\lambda_{s+1} + \beta_j + \cC_{\gP}^+) \cap \Supp M \subset (\lambda_j + \beta_j + \cC_{\gP}^+) \cap \Supp M = \emptyset$
for every $j \leq s+1$ which together with the inductive assumption completes the construction of the sequence $\lambda_1, \ldots,
\lambda_n$. The fact that $(\lambda_N + \beta_j + \cC_{\gP}^+) \cap \Supp M = \emptyset$ for every $1 \leq j \leq N$ implies that 
$(\lambda_N + \cC_{\gP}^+) \cap \Supp M = \emptyset$. Since the abelian group generated by $\Phi_\gP$ has finite index in $Q$ we 
conclude that the set $(\lambda_N + Q^+) \cap \Supp M$ is finite, where $Q^+$ is the $\Z_+$--cone generated by $\Delta^+$.
This means that the set $(\lambda_N + Q^+) \cap \Supp M$ has a maximal element $\lambda_0$ with respect to the order given
by $Q^+$, i.e. $\lambda_0 + \beta \not \in \Supp M$ for every $\beta \in \Delta^+$. 

\noindent
{\bf Case 2: $\delta \in \Delta^0$.} In this case $\Delta^\pm = (\mathring{\Delta}^\pm + \Z \delta) \cap \Delta$ and 
$\Delta^0 = (\mathring{\Delta}^0 + \Z \delta) \cap \Delta$ for some triangular decomposition $\mathring{\Delta} =
\mathring{\Delta}^- \sqcup \mathring{\Delta}^0 \sqcup  \mathring{\Delta}^+$. Notice also that $\Delta^\pm \subset \Delta^\re$. 
Since $\Delta^- \subset \Delta^i$, $\lambda \not \in \Supp M$ implies that
$(\lambda + \mathring{Q}^+ + r \Z \delta) \cap \Supp M = \emptyset$, where $\mathring{Q}^+$ is the $\Z$--cone generated by 
$\mathring{\Delta}^+$ and $\gG$ is of type $X_l^{(r)}$. It then follows that $\lambda \not \in \Supp M$ implies that
$(\lambda + (\mathring{Q}^+ \backslash \mathring{\Delta}^+) + \Z \delta) \cap \Supp M = \emptyset$. 
An argument similar to, but simpler than,
 the one used in Case 1 above shows that there exists $\lambda_0 \in \Supp M$ such that $(\lambda_0 + \mathring{Q}^+) \cap \Supp M =
 \emptyset$. This implies that $(\lambda_0 + (\mathring{Q}^+ \backslash (\mathring{\Delta}^+ + \mathring{\Delta}^+)) + \Z \delta) 
 \cap \Supp M = \emptyset$.
 
 If $\lambda_0 + \beta \not \in \Supp M$ for every $\beta \in \Delta^+$, then we are done. Otherwise,
 there exists $\lambda_1 = \lambda_0 + \beta_0 + n_0 \delta \in \Supp M$ for some $\beta_0 \in \mathring{\Delta}^+$ and some 
 $n_0 \in \Z$.
 If $\lambda_1 + \beta \not \in \Supp M$ for every $\beta \in \Delta^+$, then we are done. Otherwise,
 there exists $\lambda_2 = \lambda_1 + \beta_1 + n_1 \delta \in \Supp M$ for some $\beta_1 \in \mathring{\Delta}^+$ and some 
 $n_1 \in \Z$. Note that, for every $\beta \in \Delta^+$ we have $\lambda_2 + \beta \in 
(\lambda_0 + (\mathring{Q}^+ \backslash (\mathring{\Delta}^+ + \mathring{\Delta}^+)) + \Z \delta$, i.e. 
$\lambda_2 + \beta \not \in \Supp M$. This completes the proof.
 \hfill $\square$

\medskip

\subsection{The sets $\Delta^f$ and $\Delta^i$.} \label{subsec_sets}

\medskip

\begin{proposition} \label{proposition2.1}
Let $M \in \cW_\adm$. If $\beta_1, \beta_2 \in \Delta^f$ (respectively, $\beta_1, \beta_2 \in \Delta^i$) and $\beta \in (\Q_+ \beta_1 + \Q_+
\beta_2) \cap \Delta^\re$,
then $\beta \in \Delta^f$ (respectively, $\beta \in \Delta^i$).
\end{proposition}

\medskip

\noindent {\bf Proof.} The case $\beta_ 1 = \pm \beta_2$ is easy, so we may assume that $\beta_1$ and $\beta_2$ are not proportional.

 Consider $\Delta_0 := (\Q \beta_1 + \Q
\beta_2) \cap \Delta$ and set $\gK := \gH \oplus (\opp_{\alpha \in
\Delta_0} \gG^\alpha)$. Note that $\gK = \gK' \oplus \gZ$, where $\gK'$ is a semisimple Lie algebra of rank two
or $\gK'$ is an affine Lie algebra of type $A_1^{(1)}$ or $A_2^{(2)}$ and $\gZ$ is contained in the center of $\gK$.
Fix $\lambda \in \Supp M$ and consider $M_0
= \opp_{\mu \in Q_0} M^{\lambda + \mu}$, where $Q_0 = \langle \Delta_0 \rangle_{\Z}$ is the root lattice of $\gK$. Then $M_0 \in
\cW_\adm(\gK')$ and $\Delta^f(M_0) = \Delta^f \cap \Delta_0$,
$\Delta^i(M_0) = \Delta^i \cap \Delta_0$, where 
$\Delta^f(M_0)$ is a subset of the roots of $\gK'$ and $\Delta^f$ is a subset of the roots of $\gG$.
This implies that it suffices to
prove the theorem for the $\gK'$--module $M_0$ instead of the $\gG$--module $M$. 
For the rest of the proof we assume that $\gG$ is a semisimple Lie algebra of rank two
or $\gG$ is an affine Lie algebra of type $A_1^{(1)}$ or $A_2^{(2)}$. 

If $\beta_1, \beta_2 \in \Delta^\re$ and $\beta \in (\Q_+ \beta_1 + \Q_+ \beta_2) \cap \Delta^\re$, 
then it is obvious that $\beta_1, \beta_2 \in \Delta^i$ implies $\gamma  \in \Delta^i$.
Now assume that $\beta_1, \beta_2 \in \Delta^f$ but $\gamma \in \Delta^i$.
Note that the cone $\cC$ generated by $\Delta^i$ intersects $\Delta^\re$ in $\Delta^i$. Hence $\Delta^i = \{\pm \gamma\}$
or $\Delta^i$ is contained in the set of roots $\Delta_\gB$ of a standard Borel subalgebra of $\gG$.

Consider first the case when $\Delta^i = \{\pm \gamma\}$. Since $(\gamma, \beta_1) = (\gamma, \beta_2) = 0$ is impossible, we assume that 
$(\gamma, \beta_1) \neq 0$.  The subalgebra of $\gG$ generated by $\{e_{\beta_1}, f_{\beta_1}\}$ is isomorphic to $\gs\gl_2$, which implies that
for every $\lambda \in \Supp M$ the set 
\begin{equation} \label{eq2.1}
\left\{ \frac{2(\lambda + n \beta_1, \beta_1)}{(\beta_1, \beta_1)} \, | \, n \in \Z \right\}
\end{equation}
is a set of integers symmetric about zero. Furthermore, for any $\mu \in \gH^*$, $\mu \in \Supp M$ if and only if $\mu + \gamma \in \Supp M$.
Fix $\lambda \in \Supp M$ and consider the sets (\ref{eq2.1}) corresponding to $\lambda$ and $\lambda + \gamma$. We have
$$
\left\{ \frac{2((\lambda +\gamma) + n \beta_1, \beta_1)}{(\beta_1, \beta_1)} \, | \, n \in \Z \right\} = 
\frac{2 (\gamma, \beta_1)}{(\beta_1, \beta_1)} + \left\{ \frac{2(\lambda + n \beta_1, \beta_1)}{(\beta_1, \beta_1)} \, | \, n \in \Z \right\},
$$
which shows that they cannot both be symmetric about zero which is a contradiction with the assumption that $\Delta^i = \{\pm \gamma\}$.

Consider now the case when $\Delta^i$ is contained in the set of roots $\Delta_\gB$ of a standard Borel subalgebra of $\gG$.
By choosing appropriately the Borel subalgebra $\gB$
we can assume that
$\{\beta_1, \beta_2\}$ forms a basis of $\gB$ and 
$\gamma' = s_{\beta_1}(\beta_2) \in \Delta^i$.
Proposition \ref{prop2.0} implies that there exists $\lambda \in \Supp M$ such that $\lambda - \alpha \not \in \Supp M$ for every $\alpha \in \Delta_\gB$.
Fix a nonzero vector $v \in M^\lambda$ and denote by $M'$ the $\gG$--module generated by $v$. Then $M'$ is a highest weight module
with respect to the Borel subalgebra $\gB^{\op}$ opposite to $\gB$. Furthermore, $\Supp M' \subset \Supp M$ is finite in the direction of $\beta_1$
and $\gs \gl_2$--representation theory implies that 
$m = \frac{2(\lambda, \beta_1)}{(\beta_1, \beta_1)} \in \Z_-$
and the vector 
$v' = e_{\beta_1}^m \cdot v \in M^{s_{\beta_1}(\lambda)}$ is a highest weight vector of $M'$ with respect to the Borel subalgebra 
$s_{\beta_1}(\gB^{\op})$. Since  $\gamma' \in \Delta^i$ we conclude that 
$s_{\beta_1}(\lambda) + n \gamma' \in \Supp M'$ for every $n \in \Z_+$. Again $\gs \gl_2$--representation theory implies that $\Supp M'$
is $s_{\beta_1}$--invariant, which using that $\gamma' = s_{\beta_1}(\beta_2)$ gives
$$
s_{\beta_1} (s_{\beta_1}(\lambda) + n \gamma') = \lambda + n \beta_2 \in \Supp M' \subset \Supp M,
$$
which contradicts the assumption that $\beta_2 \in \Delta^f$. This completes the proof.
\hfill $\square$

\medskip

\begin{corollary} \label{cor2.21}
If $\gG = A_1^{(1)}$ or $\gG = A_2^{(2)}$ then $\Delta^f $ and $\Delta^i$ are one of the following
\begin{itemize}
\item[{\rm (i)}]  $\emptyset$;
\item[{\rm (ii)}]  $P \cap \Delta^{\re} $ for a principal parabolic subset $P$ of $\Delta$;
\item[{\rm(iii)}] $\Delta^{\re}$.
\end{itemize}
\end{corollary}

\medskip

\noindent
{\bf Proof.} Easy exercise. \hfill $\square$ 

\medskip

\subsection{The parabolic subalgebra $\gP_M$.} \label{subsec_decomp}

\medskip

\begin{proposition} \label{prop2.31}
Let $M \in \cW_\adm$ such that $\emptyset \neq \Delta^f$ is a proper subset of $\Delta^\re$. 
There exists a proper parabolic subset $P \subset \Delta$ such that
\begin{equation} \label{eq2.32}
P \backslash (-P) \subset \Delta^f \quad {\text { and }} (-P) \backslash P \subset \Delta^i. 
\end{equation}
\end{proposition}

\medskip

\noindent
{\bf Proof.} We consider two cases.

\noindent
{\bf Case 1. For every $\beta \in \Delta^\re$ the set $\beta + \Z \delta$ intersects both $\Delta^f$ and $\Delta^i$.}
Consider the decomposition $\mathring{\Delta} = \mathring{\Delta}' \sqcup \mathring{\Delta}''$, where
$\mathring{\Delta}' = \{ \beta \in \mathring{\Delta} \, | \, \beta + n \delta \in \Delta^f {\text { for large enough }} n \in \Z_+\}$
and 
$\mathring{\Delta}'' = \{ \beta \in \mathring{\Delta} \, | \, \beta - n \delta \in \Delta^f {\text { for large enough }} n \in \Z_+\}$.
(If $\gG$ is twisted, we assume that $\beta + n \delta \in \Delta$ in the formulas above.) Corollary \ref{cor2.21}
 implies that both sets $\mathring{\Delta}'$ and $\mathring{\Delta}''$ are
symmetric, i.e. $\mathring{\Delta}' = - \mathring{\Delta}'$ and $\mathring{\Delta}'' =- \mathring{\Delta}''$. Furthermore,
Proposition \ref{proposition2.1} implies that both $\mathring{\Delta}'$ and $\mathring{\Delta}''$ are closed. Since 
$\mathring{\Delta}$ is an irreducible root system, we conclude that $\mathring{\Delta}' = \emptyset$ or 
$\mathring{\Delta}'' = \emptyset$. Without loss of generality we assume that $\mathring{\Delta}'' = \emptyset$. 
Then we set 
\begin{equation} \label{eq2.33}
P_M = \Delta^f \cup (-\Delta^i) \cup \Z_{+} \delta.
\end{equation}
Clearly $P_M \cup (-P_M) = \Delta$. Because $\Delta^f \neq \emptyset$,  Corollary \ref{cor2.21} implies that $P_M$ is a proper subset of $\Delta$.  Furthermore, by  Proposition \ref{proposition2.1} $P_M$ is closed, and therefore a proper parabolic subset of $\Delta$. 
Finally, (\ref{eq2.32}) follows directly from the definition of 
$P_M$.

\noindent
{\bf Case 2. There exists $\beta \in \Delta^\re$ such that $\beta + n \delta$ is contained entirely in $\Delta^f$ or in $\Delta^i$.}
Set 
$$
\mathring{\Delta}^f = \{\beta \in \mathring{\Delta} \, | \, \beta + n \delta \in \Delta^f\}, \quad
\mathring{\Delta}^i = \{\beta \in \mathring{\Delta} \, | \, \beta + n \delta \in \Delta^i\}, \quad
\mathring{\Delta}^m = \mathring{\Delta} \backslash (\mathring{\Delta}^f \cup \mathring{\Delta}^i)
$$
and
$$
\mathring{P} = \mathring{\Delta}^f \cup (- \mathring{\Delta}^i) \cup \mathring{\Delta}^m.
$$
Define
\begin{equation} \label{eq2.34}
P_M = ((\mathring{P} + \Z \delta) \cap \Delta) \cup \Z \delta.
\end{equation}
Again Proposition \ref{proposition2.1} and Corollary \ref{cor2.21} imply that $\mathring{P}$ is a parabolic subset of 
$\mathring{\Delta}$. Furthermore, it is clear that 
$\mathring{P} \backslash (-\mathring{P}) \subset \mathring{\Delta}^f$ and  
$(-\mathring{P}) \backslash \mathring{P} \subset \mathring{\Delta}^i$. Then (\ref{eq2.34}) implies that 
$P_M$ is a parabolic subset of $\Delta$ which satisfies (\ref{eq2.32}). It remains to show that $P_M$ 
is a proper subset of $\Delta$. This is equivalent to showing that $\mathring{P}$ is a proper subset of $\mathring{\Delta}$.
The latter statement follows again from Proposition \ref{proposition2.1} and Corollary \ref{cor2.21}.
Indeed, the decomposition
$$
\mathring{\Delta} = (\mathring{\Delta}^f \cap (-\mathring{\Delta}^f)) \sqcup (\mathring{\Delta}^i \cap (-\mathring{\Delta}^i))
\sqcup \mathring{\Delta}^m \sqcup (\mathring{\Delta}^f \backslash (-\mathring{\Delta}^f)) \sqcup
(\mathring{\Delta}^i \backslash (-\mathring{\Delta}^i))
$$
shows that 
$$
\mathring{P} = \mathring{\Delta} \backslash (\mathring{\Delta}^i \backslash (-\mathring{\Delta}^i))
\quad {\text { and }} (\mathring{\Delta}^i \backslash (-\mathring{\Delta}^i)) = - 
(\mathring{\Delta}^f \backslash (-\mathring{\Delta}^f)).
$$
The assumption that $\mathring{P} = \mathring{\Delta}$ then leads to
$$
\mathring{\Delta} = (\mathring{\Delta}^f \cap (-\mathring{\Delta}^f)) \sqcup (\mathring{\Delta}^i \cap (-\mathring{\Delta}^i))
\sqcup \mathring{\Delta}^m.
$$
Each of the three set on the right hand side above is a symmetric closed subset of $\mathring{\Delta}$ and they commute
among themselves, i.e. the sum of two elements from different sets is never in $\mathring{\Delta}$. The fact that
$\mathring{\Delta}$ is an irreducible root system implies that two of them are empty. By assumption then we have 
that $\mathring{\Delta}$ equals either $\mathring{\Delta}^f$ or $\mathring{\Delta}^i$, which contradicts the assumption
for Case 2. \hfill $\square$

\medskip

For the rest of the paper, given a $\gG$--module $M \in \cW_\adm$, we fix the parabolic subalgebra $\gP_M$ 
corresponding to the parabolic subset $P_M$ defined by (\ref{eq2.33}) or (\ref{eq2.34}) depending to which 
of the two cases above applies.

\medskip

\noindent
{\bf Remark.} In \cite{DMP} an analogous parabolic set for finite dimensional simple Lie superalgebras 
is defined simply as $\Delta^f \cup (-\Delta^i)$. The presence of imaginary roots, however makes our situation more complicated. As a result of the classification, we will obtain that, for 
a simple module $M$, we have $P_M \cap \Delta^\re = \Delta^f \cup (-\Delta^i)$. 

\medskip

\subsection{Parabolic Induction} \label{subsec_Induction}

\medskip

\begin{theorem} \label{theorem2}

\noindent
{\rm{(i)}} Let $M$ be a simple $\gG$--module with corresponding parabolic subalgebra $\gP_M$. There exists
a simple $\gL$--module $N$, such that $M \cong V_{\gP_M} (N)$. 

\noindent
{\rm{(ii)}} If $M_1$ and $M_2$ are simple
$\gG$--modules with corresponding parabolic subalgebras $\gP_{M_1}$ and $\gP_{M_2}$, 
then $M_1 \cong M_2$ if and only if $\gP_{M_1} = \gP_{M_2}$ and $N_1 \cong N_2$.
\end{theorem}

\medskip

\noindent
{\bf Proof.} (i) If $\Delta^f = \emptyset$ or $\Delta^f = \Delta^\re$, then $\gP_M = \gG$ and there is nothing to
prove. Otherwise, Proposition \ref{prop2.31} implies that $\gP_M$ is a proper parabolic subalgebra of $\gG$,
Proposition \ref{prop2.0} shows that $M^{\gN^+} \neq 0$, and by applying Proposition \ref{ind_simple} we complete the proof.

(ii) This statement follows from the construction of the sets $P_{M_1}$ and $P_{M_2}$ and Proposition \ref{ind_simple}.
\hfill $\square$

\medskip

Theorem \ref{theorem2} reduces the problem of classifying the simple modules in $\cW_\fin$ to classifying
the ones for which $\Delta^\re = \Delta^f$ or $\Delta^\re = \Delta^i$. The former were classified in \cite{C} in 
the case when $\gG$ is untwisted and in \cite{CP2} in the case when $\gG$ is twisted. The rest of this paper 
is devoted to classifying the simple modules in $\cW_\fin$ with $\Delta^\re = \Delta^i$.

We conclude this section with some remarks. V. Futorny   studied extensively
weight modules over affine Lie algebras, see \cite{F2}, \cite{F3}, \cite{F4}, and 
\cite{F6}. In particular he has proved versions of Theorem \ref{theorem2}
under certain assumptions. He has also studied the structure of parabolically induced $\gG$--modules.
The statement of Theorem \ref{theorem2} is empty when $\gP = \gG$, i.e. when $\Delta^\re = \Delta^f$ or
$\Delta^\re = \Delta^i$. In the former case, V. Chari in \cite{C} and V. Chari and A. Pressley in \cite{CP1}, classified
all simple $\gG$--modules in $\cW_\fin$. The only case which eluded classification until now is the
case when $\Delta^\re = \Delta^i$. In Section \ref{subsec_cuspidal} below we will use Theorem \ref{reduce} to reduce
this case to the cases that are  already known. In the process we will in fact classify all simple modules
in $\cW_\fin^\Phi$ for all subsets $\Phi \subset \Delta$.

\medskip

\subsection{Length of bounded modules} \label{subsec_length}
Mathieu proved in \cite{M} that every finitely generated weight module with bounded weight multiplicities over a reductive finite
dimensional Lie algebra has finite length. This statement is no longer true for an affine algebra $\gG$ because
there are infinitely many one dimensional modules which are trivial over $[\gG, \gG]$ and on which $D$
acts as multiplication by an arbitrary scalar. 

\medskip

\begin{proposition} \label{lem_fin} 
Assume that the support of $M \in \cW_\fin$ is contained in a single $Q$--coset and there exists a constant $d$
such that $\dim M^\lambda \leq d$ for every $\lambda \in \Supp M$.
 Then
\begin{itemize}
\item[{\rm (i)}]  every finite dimensional
simple $\gG$--module is one dimensional and trivial as $[\gG,\gG]$--module;
\item[{\rm (ii)}] $M$ has finitely many infinite dimensional simple subfactors;
\item[{\rm(iii)}] $M$ has a simple submodule.
\end{itemize}
\end{proposition}

\medskip

\noindent
{\bf Proof.} (i) Let $L$ be a finite dimensional simple $\gG$--module.
Clearly $L$ is a generalized weight module, i.e.
each element of $\gH$ acts locally finitely on $L$. Since $L$ is simple, it is necessarily a weight module,
cf. \cite{DMP}. Furthermore, each $e_\alpha$ for $\alpha \in \Delta^\re$ acts nilpotently on $L$ and the 
classification of the simple integrable weight modules with finite dimensional weight spaces, \cite{C}
and \cite{CP2}, implies that $N$ is one dimensional and $[\gG, \gG]$ acts trivially on $L$.

\noindent
(ii) Let $M \in \cW_\fin$ be such that $\Supp M \subset \lambda + Q$ and assume that 
$\dim M^\mu \leq d$ for every $\mu \in \Supp M$. Every irreducible subquotient $L$ of $M$ is of one of the
following types:
\begin{itemize}
\item[(a)] finite dimensional;
\item[(b)] parabolically induced from a standard parabolic subalgebra such that 
$\Delta^i \neq \emptyset$;
\item[(c)] module with $\Delta^i = \Delta^\re$;
\item[(d)] infinite dimensional module with $\Delta^f = \Delta^\re$;
\item[(e)] infinite dimensional and parabolically induced from a non--standard parabolic subalgebra.
\end{itemize}
We will show that $M$ has finitely many irreducible
subfactors of each type (b)---(e). We discuss each of the cases above separately.

\noindent
{\bf Simple subquotients of types (b).}
Without loss of generality we may assume that $\delta \in \Delta^+$.  Assume that there are more than $d N_\gG^{l+1}$ simple
subfactors of $M$ of type (b), where $l$ is the rank of $\gG$ and $N_{\gG}$ was defined in Proposition \ref{prop_cone}. Then there are $d+1$ simple subfactors 
$L_0, \ldots, L_d$ of $M$ of type (b) with weights $\nu_i \in \Supp L_i$ for $0 \leq i \leq d$ such that 
$\nu_0 - \nu_i \in N_\gG Q$ for every $i$. Denote by $\gP_i$, $\gL_i$, and $Q_{\gP_i}$ 
the parabolic subalgebra of $\gG$, its Levi component,
and the abelian subgroup of $Q$ corresponding to $L_i$ as in Proposition \ref{prop_cone}. The fact that
$\nu_0 - \nu_i \in N_\gG Q \subset Q_{\gP_i}$ implies that, for $1 \leq i \leq d$,
\begin{equation} \label{eq2.61}
\nu_0 - \nu_i = \sum_{\beta \in \Phi_{\gP_i}} d_{i, \beta} \beta.
\end{equation}
Note that $\nu_i - \cC_{\gP_i}^+ \subset \Supp L_i$ for $0 \leq i \leq d$. Fix $s = 1 + \max\{ d_{i, \beta} \, |
\, 1 \leq i \leq d, \beta \in \Phi_{\gP_i}\}$ and consider the weight
$$
\nu = \nu_0 - s N \delta,
$$
where $N$ is the least common multiple of $|W_{\gL_1}|, \ldots, |W_{\gL_d}|$.
We claim that $\nu \in \Supp L_i$ for every $0 \leq i \leq d$. The statement is obvious for $i = 0$ and for $i>0$ we have
$$
\nu = \nu_i + (\nu_0 - \nu_i) - s N \delta = \nu_i - (s N \delta - \sum_{\beta \in \Phi_{\gP_i}} d_{i, \beta} \beta).
$$
The choice of $s$ together with equation (\ref{eq1.24}) imply that 
$$
s N \delta - \sum_{\beta \in \Phi_{\gP_i}} d_{i, \beta} \beta \in \cC_{\gP_i}^+,
$$
which shows that $\nu \in \Supp L_i$. The fact that $\nu \in \Supp L_i$ for $0 \leq i \leq d$ implies that 
$\dim M^\nu \geq d+1$ which is impossible. This contradiction shows that there are at most $d N_\gG^{l+1}$ simple
subfactors of $M$ of type (b).

\noindent
{\bf Simple subquotients of type (c).}
Since the support of every simple module of type (d) is a whole $Q$--coset, it is obvious that $M$ has at most 
$d$ simple subquotients of type (c). 

\noindent
{\bf Simple subquotients of type (d).}
The simple modules in $\cW_\fin$ with $\Delta^f = \Delta^\re$ were classified in \cite{C} and \cite{CP2}. 
Every such module is either an integrable highest weight module or a loop module. 

If $L$ is a simple 
integrable highest weight module then the central element $K$ acts on $L$ by a nonzero integer scalar 
(positive when $\delta \in \Delta^+$ and negative otherwise). On the other hand, \cite{BL} shows that 
$K$ acts trivially on every simple weight $\gG$--module with bounded weight multiplicities. Combining 
these two results we conclude that $M$ does not have any simple highest weight subquotients with
$\Delta^f = \Delta^\re$.  

In considering the integrable loop subquotients of $M$ we may restrict ourselves to the case when $\gG$ is
untwisted. Indeed, if $\gG = \cA(\gg, \sigma) = \oplus_{j \in \Z} \gg_{\bar{j}} \otimes t^j \oplus \F D \oplus \F K$
is a twisted affine algebra of type $X_n^{(r)}$
as in subsection  \ref{subsecAffine}, then its subalgebra $\gG' = \oplus_{j \in \Z} \gg_{\bar{0}} \otimes t^{rj} \oplus
\F D \oplus \F K$ is a untwisted affine algebra. Furthermore, since the root lattice $Q'$ of $\gG'$ has finite index in $Q$,
the support of $M$ is contained in the union of finitely many $Q'$--cosets. Thus it is enough to show that $M$
has finitely  many simple subquotients of type (e) in the case when $\gG$ is untwisted. For the rest of
this case we assume that $\gG$ is untwisted. The description of 
the simple integrable loop $\gG$--modules in \cite{C} implies that if $L$ is such a module, then there exists 
a Borel subalgebra $\gb$ of $\gg$ and a $\gb$--dominant weight
$\lambda \in \gh^*$, and a constant $a \in \F$ such that 
$$
(\Supp V(\lambda) \backslash \mathring{W}(\lambda)) + (\Z + a) \delta \subset \Supp L,
$$
where $V(\lambda)$ denotes the finite dimensional $\gg$--module with highest weight $\lambda$.
Moreover, the set $\Supp V(\lambda) \backslash \mathring{W}(\lambda)$ is not empty and thus
it contains a miniscule weight of $\gg$. Since there are only finitely many miniscule weights of $\gg$ we conclude
that there are finitely many simple integrable loop subfactors of $M$.

\noindent
{\bf Simple subqotients of type (e).} As in the previous case we assume that $\gG$ is untwisted. Consider the simple subquotients of $M$ which are parabolically induced from non--standard parabolic subalgebras 
corresponding to a fixed triangular decomposition of $\Delta$. Since $\delta \in \Delta^0$, the triangular decomposition of
$\Delta$ induces a triangular decomposition of $\mathring{\Delta}$. If $L$ is a simple subquotient of $M$
parabolically induced from a parabolic corresponding to the fixed triangular decomposition of $\Delta$,
and $\mu \in \Supp L$, then
$$
\mu - \langle \mathring{\Delta}^{+} \rangle_{\Z_+} + \Z \delta \subset \Supp L.
$$
This equation shows that there are at most $d$ simple subquotients of $M$ corresponding to parabolic
subalgebras with the fixed triangular decomposition. Since there are only finitely many triangular decompositions
of $\Delta$ corresponding to non--standard parabolic subalgebras, we conclude that $M$ has finitely many 
simple subquotients of type (e).

This completes the proof of (ii).

\noindent
(iii) Since $M$ has only finitely many infinite dimensional simple subquotients, we see that either $M$ has 
a simple infinite dimensional submodule or a submodule $M'$ each simple subfactor of which is one dimensional
and, by (i), trivial as $[\gG,\gG]$--module. In the latter case it is obvious that $M'$ is completely reducible. Indeed, since
 $M$ is a weight module all $e_{\alpha}$ and $f_{\alpha}$ act trivially on $M$. 
 Then, since $D$ acts diagonally on $\gG$, $\gG$ decomposes as a direct sum of the eigenspaces of $D$, which
implies that $M'$ (and consequently, $M$) has a one dimensional simple submodule.
\hfill $\square$

\medskip

\section{Localization}

\medskip

\subsection{Localization and parabolic induction}
In this subsection we will show that the twisted localization
correspondence $\cD_{\alpha}^x$ commutes with the parabolic
induction functors $M_{\gP}$ and $V_{\gP}$.  

Let $\gP$ be a parabolic subalgebra of $\gG$ corresponding to the
parabolic subset $P$ of $\Delta$. As usual $\gL$, $\gN^+$, and $\gN^-$ denote
the Levi component, the nillradical of $\gP$, and the nillradical of the parabolic subalgebra $\gP^{\op}$
opposite to $\gP$. Let $\Delta_\gL = P \cap (-P)$
be the roots of $\gL$ and let $Q_\gL = \langle \Delta_\gL \rangle_\Z$ be  the
root lattice of $\gL$. Let $S$ be a weight $\gL$--module whose
support is contained in a single $Q_\gL$--coset $\lambda + Q_\gL \in
\gH^*/Q_\gL$. We consider $S$ as a $\gP$--module
 with trivial action of $\gN^{+}$. 

\medskip 

\begin{proposition} \label{loc_ind}
Let $\alpha \in \Delta_\gL$ and let $S$ be an
$f_{\alpha}$--injective weight module whose support is included in
a single $Q_\gL$-coset. Then for every $x$ in $\F$

\noindent \rm{(i)} $f_{\alpha}$ acts injectively on $M_{\gP}(S)$
and $V_{\gP}(S)$;

\noindent \rm{(ii)} $\cD_{\alpha}^{x}M_{\gP}(S) \simeq
M_{\gP}(\cD_{\alpha}^{x} S)$;

\noindent \rm{(iii)} $\cD_{\alpha}^{x}V_{\gP}(S) \simeq
V_{\gP}(\cD_{\alpha}^{x} S)$.

\end{proposition}

\medskip

\noindent {\bf Proof.} We follow the same reasoning as the one in
the proof of Lemma 2.4 in \cite{G}. For completeness
we provide a sketch of the proof.

\rm{(i)} The injectivity of $f_{\alpha}$ on $M_{\gP}(S)$ follows from 
 $[f_{\alpha}, \gN^{-}] \subseteq \gN^{-}$ and the isomorphism 
 $M_{\gP}(S) \simeq U(\gN^{-}) \otimes S$ of $\gN^-$--modules.

To prove that $f_\alpha$ acts injectively on $V_\gP(S)$, assume the 
contrary. 
Then there exists $m \in M_{\gP}(S)$ such that 
$m \notin Z_{\gP}(S)$ and $f_{\alpha}\cdot m \in Z_{\gP}(S)$. 
Since $Z_{\gP}(S)$ is
$f_{\alpha}$--injective, we have  $m \in \cD_{\alpha}Z_{S} \cap M_{\gP}(S) \subset \cD_{\alpha} M_{\gP}(S)$.
On the other
hand, the maximality of $Z_{\gP}(S)$ and the inclusion $Z_{\gP}(S) \subseteq
\cD_{\alpha}Z_{\gP}(S) \cap M_{\gP}(S) \subset M_{\gP}(S)$ imply
that $Z_{\gP}(S) = \cD_{\alpha}Z_{\gP}(S) \cap M_{\gP}(S)$. Thus
$m \in Z_{\gP}(S)$ which is a contradiction.

\rm{(ii)} We use again that $[f_{\alpha}, \gN^{-}] \subseteq
\gN^{-}$. The statement follows by verifying that
\begin{equation} \label{loc_n}
\cD_{\alpha}(U(\gN^{-})\cdot S) = U(\gN^{-}) \cdot \cD_{\alpha}(S)
\end{equation}
and
\begin{equation}  \label{phi_n}
\Phi_{\alpha}^x (U(\gN^{-})\cdot R) = U(\gN^{-}) \cdot
\Phi_{\alpha}^x (R)
\end{equation}
for any $U_{\alpha}$--module $R$ (and in particular for
$R=\cD_{\alpha}S$). The two identities above should be understood
as identities of $\gG^{\alpha}$--modules (see Remark
\ref{rem_loc}). To prove the identity (\ref{loc_n}) we apply
several times the formula
$$
f_{\alpha}^n u f_{\alpha}^{-n}= \sum\limits_{i=0}^\infty
\binom{n}{i}\, \ad(f_\alpha)^i (u) \, f_\alpha^{-i}
$$
for $n \in \Z$. For (\ref{phi_n}) we use (\ref{theta}).

\rm{(iii)} We first show that $\cD_{\alpha}^x Z_{\gP}(S)$ is a
maximal submodule of $\cD_{\alpha}^x M_{\gP}(S) \simeq M_{\gP}
(\cD_{\alpha}^x S)$ which has a trivial intersection with $\cD_{\alpha}^x
S$. Indeed, suppose that there is a $U_{\alpha}$-module $Z'$ for
which
$$ \cD_{\alpha}^x Z_{\gP}(S)
\subseteq Z' \subset \cD_{\alpha}^x M_{\gP}(S).
$$
Then by applying $\Phi_{\alpha}^{-x}$ to the inclusions above we
obtain
$$
\cD_{\alpha} Z_{\gP}(S) \subseteq \Phi_{\alpha}^{-x} Z' \subset
\cD_{\alpha} M_{\gP} (S).
$$
Now intersecting with $M_{\gP}(S)$ and using that $Z_{\gp}(S) =
\cD_{\alpha}Z_{\gP}(S) \cap M_{\gP}(S)$ (proved in (i)) we
conclude that
$$
Z_{\gp}(S) \subseteq \Phi_{\alpha}^{-x} Z' \cap M_{\gP} (S)
\subset M_{\gP}(S)
$$
By the maximality of $Z_{\gP}(S)$ we have $Z_{\gp}(S) =
\Phi_{\alpha}^{-x} Z' \cap M_{\gP} (S)$. This easily implies that
$\cD_{\alpha} Z_{\gp}(S) = \Phi_{\alpha}^{-x} Z'$ (and hence
$Z'=\cD_{\alpha}^x Z_{\gp}(S)$). Indeed, if $z \in
\Phi_{\alpha}^{-x} Z'$ then for some $N$, $f_{\alpha}^Nz \in
M_{\gP}(S) \cap \Phi_{\alpha}^{-x} Z' = Z_{\gP}(S)$, and thus $z
\in \cD_{\alpha}Z_{\gP}(S)$. It remains to show that
$$
\cD_{\alpha}^x \left(  M_{\gP}(S)/ Z_{\gP}(S)\right) \simeq
\cD_{\alpha}^x \left(  M_{\gP}(S)\right) / \cD_{\alpha}^x \left(
Z_{\gP}(S)\right).
$$
This follows from the fact that $M_{\gP}(S)/ Z_{\gP}(S)$ is
$f_{\alpha}$--injective (see Lemma 2.2 (ii) in \cite{G} for the
complete proof).
 \hfill $\square$

\medskip

\subsection{Localization and loop modules} \label{subsec2.31}

\medskip

\begin{proposition} \label{loc_tensor}
Let $\gG = \cL(\gg)$ be a untwisted affine Lie algebra,
let $\alpha \in \mathring{\Delta}$, and let 
${\cL}_{a_0,...,a_k}(L_0 \otimes F_1 \otimes...\otimes F_k)$ be a loop
$\gG$-module for which $a_i \in \C$, $F_i$ are simple
finite dimensional $\gg$-modules, and $L_0$ is a simple
$f_{\alpha}$-injective $\gg$-module. Then
$$
{\cD}_{\alpha + r \delta}{\cL}_{a_0,...,a_k}(L_0
\otimes F_1 \otimes...\otimes F_k) \simeq {\cL}_{a_0,...,a_k}({\cD}_{\alpha}L_0 \otimes F_1
\otimes...\otimes F_k)
$$
\end{proposition}

\medskip

\noindent {\bf Proof.} 
The proof follows the reasoning of the proof of Theorem 3.1 in \cite{De}. For the sake of completeness we outline the important steps. It can be easily verified that ${\cD}_{\alpha + r \delta}{\cL}_{a_0,...,a_k}(L_0 \otimes F_1 \otimes...\otimes F_k)$ is $f_{\alpha + r \delta}$-injective module for every integer $r$.  We define the map 
$$
g: {\cD}_{\alpha + r \delta}{\cL}_{a_0,...,a_k}(L_0
\otimes F_1 \otimes...\otimes F_k) \to {\cL}_{a_0,...,a_k}({\cD}_{\alpha}L_0 \otimes F_1
\otimes...\otimes F_k)
$$ 
by
\begin{equation} \label{loc_loop_iso}
\left( f_{\alpha} \otimes t^r \right)^{-N} (v_0 \otimes ... \otimes v_k \otimes t^s) \mapsto \sum_{\substack{i_1,...,i_k\geq 0 \\ i_0 = -N - i_1-...-i_k}} \binom{-N}{i_1,...,i_k} a_0^{i_0 r} f_{\alpha}^{i_0}v_0 \otimes ... \otimes a_k^{i_kr} f_{\alpha}^{i_k}v_k \otimes t^{-Nr+s},
\end{equation}
where $\binom{n}{k_1,...,k_j}:=\binom{n}{k_1}\binom{n-k_1}{k_2}...\binom{n-k_1-...-k_{j-1}}{k_j}$ for $n \in Z$ and $k_r \in \Z_{+}$. In fact,  $\binom{n}{k_1,...,k_j}$ is exactly the coefficient of $a_0^{n-k_1-...-k_j} a_1^{k_1}...a_j^{k_j}$ in the expansion of $(a_0 +...+ a_j)^n$.

We first notice that $g = \Id$ for $N \leq 0$ as it is easy to check that for $n \in \Z_{+}$
$$
(f_{\alpha} \otimes t^r)^n  (v_0 \otimes...\otimes v_k
\otimes t^s)
=\sum_{\substack{ i_0,...,i_k \geq 0 \\ i_0+...+i_k = n}}\frac{n!}{i_0!...i_k!}a_0^{i_0r} f_{\alpha}^{i_0}v_0
\otimes ... \otimes a_k^{i_kr}f_{\alpha}^{i_k}v_k  \otimes t^{nr+s}.
$$
Next we show that $g$ is well-defined. Let 
$$
\left( f_{\alpha} \otimes t^r \right)^{-N} (\sum_{i} v_0^i \otimes ... \otimes v_k^i \otimes t^{s_i}) =  
\left( f_{\alpha} \otimes t^r \right)^{-K} (\sum_{j} w_0^j \otimes ... \otimes w_k^j \otimes t^{r_j})
$$ 
for $N,K \geq 0$. Then 
$$
\left( f_{\alpha} \otimes t^r \right)^{K} (\sum_{i} v_0^i \otimes ... \otimes v_k^i \otimes t^{s_i}) =  
\left( f_{\alpha} \otimes t^r \right)^{N} (\sum_{j} w_0^j \otimes ... \otimes w_k^j \otimes t^{r_j}).
$$
Furthermore, we have that
\begin{align*}
  &  \left( f_{\alpha} \otimes t^r \right)^{N+K} g   \left( f_{\alpha} \otimes t^r \right)^{-N} \left( v_0 \otimes...\otimes v_k \otimes t^s \right) \\
\;  &=  \left( f_{\alpha} \otimes t^r \right)^{N+K} \sum_{\substack{i_1,...,i_k\geq 0 \\ i_0 = -N - i_1-...-i_k}} \binom{-N}{i_1,...,i_k} a_0^{i_0 r} f_{\alpha}^{i_0}v_0 \otimes ... \otimes a_k^{i_k r} f_{\alpha}^{i_k}v_k \otimes t^{-Nr+s} \\
  &= \sum_{\substack{i_1,...,i_k\geq 0 \\ i_0 = -N - i_1-...-i_k}} \sum_{\substack{j_1,...,j_k\geq 0 \\ j_0 +...+ j_k = N + K}} \binom{-N}{i_1,...,i_k} \binom{N + K}{j_1,...,j_k} a_0^{(i_0+j_0) r} f_{\alpha}^{i_0+ j_0}v_0 \otimes ... \otimes a_k^{(i_k + j_k)r} f_{\alpha}^{i_k+j_k}v_k \otimes t^{Kr+s} \\
& =  \sum_{\substack{l_1,...,l_k\geq 0 \\ l_0 = -N - l_1-...-l_k}}  \binom{-K}{l_1,...,l_k} a_0^{l_0 r} f_{\alpha}^{l_0}v_0 \otimes ... \otimes a_k^{l_k r} f_{\alpha}^{l_k }v_k \otimes t^{Kr+s} \\
&=  \left( f_{\alpha} \otimes t^r \right)^{K} \left( v_0 \otimes...\otimes v_k \otimes t^s \right) .
\end{align*}
(Here we  used the fact that 
$$\sum_{\substack{ i_1,..,i_k, j_1,...,j_k \geq 0 \\ i_t+j_t = l_t} }\binom{-N}{i_1,...,i_k} \binom{N + K}{j_1,...,j_k} = \binom{K}{l_1,...,l_k},$$ which follows by comparing the coefficients of $a_0^{l_0}...a_k^{l_k}$ in the identity $(a_0+...+a_k)^{-N} (a_0+...+a_k)^{N+K } = (a_0+...+a_k)^{K}$.)
Using the identities above we conclude that
\begin{align*}
& \left( f_{\alpha} \otimes t^r \right)^{N+K} g\left( \left( f_{\alpha} \otimes t^r \right)^{-N}  \left( \sum_{i} v_0^i \otimes ... \otimes v_k^i \otimes t^{s_i} \right) \right) \\
& = \left( f_{\alpha} \otimes t^r \right)^{K} \sum_i v_0^i \otimes...\otimes v_k^i \otimes t^s_i\\
& = \left( f_{\alpha} \otimes t^r \right)^{N} (\sum_{j} w_0^j \otimes ... \otimes w_k^j \otimes t^{r_j}) \\
& = \left( f_{\alpha} \otimes t^r \right)^{N+K} g\left( \left( f_{\alpha} \otimes t^r \right)^{-K} \left( \sum_{j} w_0^j \otimes ... \otimes w_k^j \otimes t^{r_j} \right) \right).
\end{align*}
This completes the proof that $g$ is well defined. Clearly, $g$ is a homomorphism of $\cL(\gg)$--modules.

Finally, to show that $g$ is an isomorphism of $\cL (\gg)$-modules we construct the map inverse to $g$. It is defined by
\begin{align*}
\tilde{g}:  & f_{\alpha}^{-n} v_0 \otimes v_1 \otimes...\otimes v_k \otimes t^s \mapsto  \\ & \left( f_{\alpha} \otimes t^r \right)^{-N}\left( \sum_{\substack{i_1,...,i_k\geq 0 \\ i_0+ i_1+ ... + i_k = N}} \binom{N}{i_1,...,i_k} a_0^{i_0 r} f_{\alpha}^{i_0 - n}v_0 \otimes ... \otimes a_k^{i_kr} f_{\alpha}^{i_k}v_k \otimes t^{-Nr+s} \right),
\end{align*}
where $N$ is chosen as follows. Let $m_i$ be such that $f_{\alpha}^{m_i}$ act trivially on $F_i$. Then choose any $N$ with $N \geq m_1 +...+ m_k +n$. The choice of $N$ leads to $i_0 - n = N-i_1-...-i_k - n \geq 0$ in the summation above. We leave it to the reader to verify that $\tilde{g}$ is indeed inverse to 
$g$. \hfill $\square$

\medskip

\subsection{Main theorem on localization}\label{subsec_length}

\begin{proposition} \label{theorem21}
Let $M$ be a simple module in ${\mathcal W}_{fin}$ and let $\alpha
\in \Delta^\re$ be such that $f_{\alpha}$ and $e_{\alpha}$ act
injectively on $M$. If $e_\alpha$ acts nilpotently on at least one
element of $M' = \Phi_\alpha^x M$, then $M'$ has a simple
submodule on which $e_\alpha$ acts locally nilpotently.
\end{proposition}

\medskip

\noindent {\bf Proof.}
First we consider the case when $\Delta_M^f = \emptyset$. Then all weight spaces
of $M$ has the same finite dimension and so do the weight spaces of $\Phi_\alpha^x M$.
Fix a nonzero vector $m \in \Phi_\alpha^x M$ on which $e_\alpha$ acts nilpotently
and denote by $M''$ the submodule of $\Phi_\alpha^x M$ generated by $m$. Then
$M''$ has  bounded weight multiplicities and $e_\alpha$ acts locally nilpotently on
$M''$. Applying Proposition \ref{lem_fin} we conclude that there exists a simple submodule $M'$ of
$M''$. Clearly $M'$ is a submodule of $\Phi_\alpha^x M$ with the desired properties.

Assume now that
$\Delta_M^f \neq \emptyset$. Fix a triangular decomposition
$\Delta = \Delta^{+} \sqcup \Delta^0 \sqcup \Delta^{-}$ and a
simple weight $\gL$-module $S$ such that $M \simeq L_{\gP}(S)$.
We have two possibilities for $\delta$ with respect to $\Delta^0$.

\noindent
{\bf Case 1: $\delta \notin \Delta^0$.} Now $\gL$ is a
finite-dimensional reductive Lie algebra. Then, $\Phi_{\alpha}^x
S$ as bounded $\gL$-module has finite length. Let $S'$ be a simple
submodule. We know that ${\cD}_{\alpha}^{-x}S'  \simeq S$. Then
by Proposition \ref{loc_ind}, ${\cD}_{\alpha}^{-x}V_{\gP}(S')
\simeq V_{\gP}(S) \simeq M$ and therefore $\Phi_{\alpha}^x M
\simeq {\cD}_{\alpha}V_{\gP}(S') $ contains a simple submodule
isomorphic to $V_{\gP}(S')$.

\noindent
{\bf Case 2: $\delta \in \Delta^0$} In this case $\gL = \gL_1
\oplus ... \oplus \gL_r$. Then $S = S_1 \otimes ... \otimes S_r$.
Assume that $\alpha \in \Delta_{\gL_1}$. Then $\Phi_{\alpha}^x S \simeq
(\Phi_{\alpha}^x S_1) \otimes... \otimes S_r$.
Since $\Delta_{S_1}^I = \Delta_{S_1}^{\re}$, from Lemma
\ref{lem_fin} and Case 3.1 we conclude that there is a simple
submodule $S_1'$ of $\Phi_{\alpha}^x S_1$. In particular, ${\cD}_{\alpha}^{-x} S_1' \simeq S_1$. Then 
${\cD}_{\alpha}^{-x}(S_1' \otimes ... \otimes S_r) \simeq ({\cD}_{\alpha}^{-x} S_1') \otimes ... \otimes S_r = S$. As in Case
3.1 we conclude that $V_{\gP}(S_1' \otimes ... \otimes S_r)$ is a
simple submodule of $\Phi_{\alpha}^x M \simeq {\cD}_{\alpha}
V_{\gP}(S_1' \otimes ... \otimes S_r)$.
\hfill $\square$

\medskip

We are now ready to prove the main Theorem of this section.

\medskip

\begin{theorem}\label{reduce}
Let $M$ be a simple module in ${\mathcal W}_{\fin}$ and let $\pm\alpha \in \Delta^i(M)$.
There exist $x \in \F$ and simple modules $N, L \in \cW_\fin$ such that
$M \simeq \cD_\alpha^x N \simeq \cD_{-\alpha}^x L$ and $\alpha \in \Delta^f(N)$ and $-\alpha \in \Delta^f(L)$.
Moreover, if $\beta \in \Delta^f(M)$ then $\beta \in \Delta^f(N)$ and $\beta \in \Delta^f(L)$.
\end{theorem}

\noindent
{\bf Proof.} Let $\lambda \in \Supp M$ and let $v \in M^{\lambda}$ be an eigenvector of $f_{\alpha}e_{\alpha}$, i.e.
$f_{\alpha}e_{\alpha}v = c v$ for some $c \in \F$. We fix $x \in \F$ to be one of the roots of the quadratic
equation $c - z\langle\lambda, \alpha\rangle + z^2\langle \alpha, \alpha\rangle = 0$. Set $N':= \cD_\alpha^x M$. Then
$$
e_{\alpha}(f_{\alpha}^{x}\cdot v)= f_{\alpha}^{x}\cdot
(\Theta_x(e_{\alpha}) \cdot v).
$$
On the other hand
$$
\Theta_x(e_{\alpha})  = \sum_{i\geq 0}
\binom{x}{i}\,\mbox{ad}(f_{\alpha})^{i}(e_{\alpha})f_{\alpha}^{-i}
= e_{\alpha} - xh_{\alpha}f_{\alpha}^{-1} +\binom{x}{2}
[f_{\alpha}, h_{\alpha}]f_{\alpha}^{-1}.
$$
Since $e_{\alpha}v = c f_{\alpha}^{-1}v$, we conclude
$$
e_{\alpha}(f_{\alpha}^{x}\cdot v) = f_{\alpha}^{x}\cdot
((e_{\alpha} - x \langle\lambda - \alpha,\alpha \rangle
f_{\alpha}^{-1}+x(x-1)\langle\alpha,\alpha\rangle
f_{\alpha}^{-1})v) = 0.
$$
Set $N'':=\{n \in N' \; | \; e_{\alpha}^a n = 0 \mbox { for some }
a>0\} $.  Then $N''$ is a nontrivial submodule of $N'$. This
implies that  $\cD_\alpha^{-x} N' \subseteq \cD_\alpha^{-x} N''
\subseteq M$ and since $M$ is simple we have $\cD_\alpha^{-x} N' =
\cD_\alpha^{-x} N'' = M$. By Proposition \ref{theorem21} we know
that $N'$ has a simple nontrivial submodule $N$. Using that $M$ is
simple we conclude that $\cD_\alpha^x N = M$.

For the last part of the proposition, assume that $\beta \in \Delta^f (M)$ and consider $\mu \in \Supp N$. Then
$\Phi_\alpha^{-x} (\cD_\alpha N) = \cD_\alpha^{-x} N = M$ implies  that $\mu - x\alpha \in \Supp M$.
Proposition \ref{prop3} implies that there is $n_0 \geq 0$ such that $\mu - x\alpha + n \alpha \notin \Supp M$ for
every $n\geq n_0$. Therefore $\mu + n \alpha \notin \Supp \cD_\alpha^x M$ for $n \geq n_0$, i.e. $\alpha \in \Delta^f(N)$.
\hfill $\square$

\medskip

\begin{proposition} \label{prop4}
Let $M$ and $N$ be two simple modules on both of which $f_\alpha$ acts injectively. Then
$\cD_\alpha M \simeq \cD_\alpha N$ implies $M \simeq N$.
\end{proposition}

\noindent
{\bf Proof.} We may consider $M$ and $N$ as submodules
of a module isomorphic to $\cD_\alpha M \simeq \cD_\alpha N$.  Fix
$0 \neq v \in M$. There exists $n \in \Z_+$ such that $f_\alpha^n v \in N$. Then $0 \neq f_\alpha^n v \in M \cap N$ and since
both $M$ and $N$ are simple, we conclude that $M=N$.
\hfill $\square$

\medskip

\subsection{Localization of imaginary Verma $\widehat{\gs \gl (2)}$-modules}

\medskip

\begin{lemma} \label{e-f-loc}
Let $M$ be a simple module in ${\mathcal W}_{fin}$ and $\alpha \in \Delta^{\re}$ be such that $f_{\alpha}$ acts injectively on $M$.
Let also $v  \in M^{\lambda}$ be such that $e_{\alpha} v = 0$. Then there is $x \in \F$ such that 
$ e_{\alpha}^k(f_{\alpha}^x \cdot v) = c_{k,x} f_{\alpha}^{x-k} \cdot v$  for every $k \in \Z_{+}$ and the constants $c_{k,x}$ are such that $c_{k, x+l} \neq 0$ for every $k \in \Z_{+}$ and $l \in \Z$.
\end{lemma}

\medskip

\noindent {\bf Proof.} Using Lemma \ref{lmnew} we verify that
\begin{eqnarray*}
e_{\alpha} (f_{\alpha}^x \cdot v) & = & f_{\alpha}^x \cdot (\Theta_x (e_{\alpha} )v)\\
 & = & f_{\alpha}^x \cdot \left(  \sum\limits_{i=0}^\infty \binom{x}{i}\, \ad(f_\alpha)^i (e_{\alpha}) \, f_\alpha^{-i}v\right) \\
 & = & f_{\alpha}^x \cdot \left(e_{\alpha}v - x h_{\alpha} f_{\alpha}^{-1}v  - \binom{x}{2}\langle \alpha, \alpha \rangle f_{\alpha}^{-1} v\right) \\
 & = & \left(- x \langle \lambda - \alpha, \alpha \rangle   - \binom{x}{2}\langle \alpha, \alpha \rangle \right) f_{\alpha}^x \cdot \left(f_{\alpha}^{-1} v \right) \\
 & = & p_{\lambda}(x) f_{\alpha}^{x-1} \cdot  v,
\end{eqnarray*}
where $p_{\lambda}(x)= - x \langle \lambda - \alpha, \alpha \rangle   - \binom{x}{2}\langle \alpha, \alpha \rangle$. Now we easily conclude that
$$
e_{\alpha}^k \left( f_{\alpha}^x \cdot v\right) = p_{\lambda}(x) p_{\lambda}(x-1)... p_{\lambda}(x-k+1)f_{\alpha}^{x-k} \cdot  v 
$$
To complete the solutions we choose any $x$ that is not contained in $\{x_1 + \Z \}\cup \{x_2 + \Z \} $, where $x_1$ and $x_2$ are the roots of $p_{\lambda}(x)=0$.
 \hfill $\square$
 
 \medskip

For the rest of this section we fix $\gg = \gs \gl (2)$, a real root $\alpha$, and a standard $\gs \gl (2)$-triple $\{e, f, h\}$ 
for which $e_{\alpha} = e\otimes t^0$. Consider the imaginary triangular decomposition 
${\mathcal N}^+ \oplus {\mathcal H} \oplus {\mathcal N}^-$ of ${\mathcal L} (\gg)$. Namely
${\mathcal N}^{+}$ is spanned by $e \otimes t^n$, ${\mathcal N}^{-}$  by $f \otimes t^n$,  and $\mathcal H$ by $h \otimes t^n$. 
The imaginary Verma ${\mathcal L} (\gg)$-module is defined as 
$M(\lambda) = U({\mathcal L}(\gg)) \otimes_{U({\mathcal H} \oplus {\mathcal N}^+)} \F v $, where $v$ is such that 
${\mathcal N}^+\cdot v = 0$, $(h\otimes t^n)\cdot v = 0$ for $n\neq 0$ and $(h\otimes t^0)\cdot v = \lambda v$. 
A result of Futorny (see \cite{F2}) implies that $M(\lambda)$ is simple if $\lambda \neq 0$. The next proposition 
shows that if we drop the finiteness condition for the weight multiplicities, we  may have a simple weight 
$\gG$-module that remain simple if considered as a ${\mathcal L}(\gg)$-module.

\medskip

\begin{proposition} \label{rao}
Let $x \in \F$ be such that Lemma \ref{e-f-loc}  holds. Then  ${\mathcal D}_x^{\alpha}M(\lambda)$ is a simple 
${\mathcal L}(\gg)$-module.
\end{proposition}

\medskip

\noindent {\bf Proof.} Put for simplicity $U:= U({\mathcal L}(\gg))$, $f_n:= f \otimes t^n$, and  $e_n:= e \otimes t^n$. In particular, $e_{\alpha} = e_0$  and $f_{\alpha} = f_0$. We prove the theorem in two steps.

{\it Step 1: ${\mathcal D}_x^{\alpha}M(\lambda)$ is generated by $f^x \cdot v$, i.e. ${\mathcal D}_x^{\alpha}M(\lambda) = U \cdot (f^x \cdot v)$}.
We first fix the following basis of $M(\lambda)$
$${\mathcal B} = \{ f_{k_1}^{n_1}f_{k_2}^{n_2}...f_{k_r}^{n_r}\cdot v \; | \; k_i \in \Z, k_1< ...< k_r, n_i \in \Z_{+}\}.$$
Using the basis ${\mathcal B}$, for any $x \in \F$, we easily find bases ${\mathcal B}_{\alpha}$ and  ${\mathcal B}_{\alpha}^x$  for ${\mathcal D}_{\alpha} M(\lambda)$ and ${\mathcal D}_{\alpha}^x M(\lambda)$, respectively. Namely:
$${\mathcal B}_{\alpha} = \{ f_0^{n_0} f_{k_1}^{n_1}...f_{k_r}^{n_r} \cdot v \; | \; k_i, n_0 \in \Z, k_1< ...< k_r,k_i \neq 0, n_l\in \Z_{+}\mbox{ for }l>0 \}$$
and
$${\mathcal B}_{\alpha}^x =  \{  f_0^{n_0} f_{k_1}^{n_1}...f_{k_r}^{n_r} \cdot v \; | \; k_i, n_0 \in \Z, k_1< ...< k_r, k_i \neq 0,n_l\in \Z_{+}\mbox{ for }l>0\}.$$
It is thus enough to show that $ f_0^{n_0} f_{k_1}^{n_1}...f_{k_r}^{n_r} \cdot v$ is in $U \cdot (f^x \cdot v)$  for $n_0 \in \Z$, $x \in \F$ any positive integers $n_1,...,n_r$, and any integers $k_1,...,k_r$. This is obvious for $n_0 \geq 0$. In the case of negative $n_0$ we use that 
$$f_0^{n_0+x} f_{k_1}^{n_1}...f_{k_r}^{n_r} \cdot v  = f_0^{n_0} f_{k_1}^{n_1}...f_{k_r}^{n_r} \cdot (f_0^x \cdot v).  $$
Now applying Lemma \ref{e-f-loc} we conclude that $ f_0^{n_0+x} f_{k_1}^{n_1}...f_{k_r}^{n_r}  \cdot v  = u \cdot (f_0^x \cdot v)$ for $u = \frac{1}{c_{-n_0, x}}e_0 ^{-n_0} f_{k_1}^{n_1}...f_{k_r}^{n_r}  \in U$.

{\it Step 2: For any element $m^x$ in ${\mathcal D}_x^{\alpha}M(\lambda)$ there is $u$ in $U$ for which $u \cdot m^x = f_0^x \cdot v$.} We first notice that any element $m^x$ in ${\mathcal D}_x^{\alpha}M(\lambda)$ has the form $m^x = f_0^{x-k}\cdot m$ for some $m \in M(\lambda)$  and nonnegative $k$. Then we fix an element $u^-$ in $U({\mathcal N}^-)$ for which  $u^-\cdot m = v$.  Now using Lemma \ref{lmnew} we find
\begin{eqnarray*}
f_0^x \cdot v  & = & f_0^x \cdot (u^- \cdot m)\\
 & = & f_0^x u^- f_0^{-x} \cdot (f_0^x \cdot m)\\
 & = & \Theta_{-x}(u^-) \cdot f_0^k  \cdot (f_0^{x-k} \cdot m).
\end{eqnarray*}
Let us fix positive integer $l$ for which $\Theta_{-x}(u^-) = f_0^{-l}u_1$ for some $u_1 \in U$. Then by Lemma \ref{e-f-loc} we have that 
$$
f_0^x \cdot v = \frac{1}{c_{l,x+l}}e_0^l \cdot (f_0^{x+l} \cdot v).
$$
Then we easily verify that for $u = \frac{1}{c_{l,x+l}}e_0^l u_1 f_0^k \in U$ we have 
$$u \cdot (f_0^{x-k}\cdot m) = \frac{1}{c_{l,x+l}}e_0^l (u_1 \cdot f_0^{x} \cdot v) = \frac{1}{c_{l,x+l}}e_0^l \cdot f_0^{x+l}\cdot v = f_0^x \cdot v$$
\hfill $\square$

\medskip

\section{Simple $\gG$--modules with finite dimensional weight spaces} \label{lists}
To complete the classification we need to describe which parabolically induced modules have finite dimensional weight spaces
and, in view of Theorem \ref{reduce}, to determine which parabolically induced modules have bounded weight multiplicities and to
determine the cuspidal modules we obtain from them by twisted localization. 

\medskip

\subsection{Standard and mixed type parabolic subalgebras} \label{standard_and_mixed}

\medskip
First we show that mixed type parabolic subalgebras do not yield any simple weight modules with finite dimensional weight spaces.

\begin{proposition} \label{prop4.1}
Let $\gP$ be a mixed type parabolic subalgebra of $\gG$ with Levi component $\gL$
and let $N$ be a simple weight $\gL$--module with finite dimensional weight spaces.
Then $V_\gP(N) \not \in \cW_\fin$. 
\end{proposition}

\medskip

\noindent
{\bf Proof.} $\gG$ contains a subalgebra $\gK$ isomorphic to $A_1^{(1)}$ such that $\gP \cap \gK$ is a mixed type Borel subalgebra of $\gK$.
Consider a non--trivial $\gK$--submodule $M'$ of $V_\gP(N)$. V. Futorny, \cite{F2}, established that every non--trivial simple
$\gP \cap \gK$--highest weight module has (some) infinite dimensional weight spaces, which implies that $V_\gP(N) \not \in \cW_\fin$.
\hfill $\square$

\medskip

Our next step is to show that non trivial parabolically induced modules from standard parabolic subalgebras do not have bounded weight multiplicities.

\begin{proposition} \label{prop4.2}
Let $\gP$ be a standard parabolic subalgebra of $\gG$ with Levi component $\gL$
and let $N$ be a simple weight $\gL$--module with finite dimensional weight spaces.
Then the weight multiplicities of $V_\gP(N)$ are not bounded.
\end{proposition}

\medskip

\noindent
{\bf Proof.} If $K$ acts nontrivially on $N$, then the statement follows from the result of Britten and Lemire, \cite{BL}.
To complete the proof it is enough to prove the statement for a highest weight module of $A_1^{(1)}$ on which $K$ acts trivially.

Let $\gG \cong A_1^{(1)}$ and fix a root basis $e_n = e \otimes t^n, f = f \otimes t^n, h_n = h \otimes t^n, D, K$ of $\gG$, where
$e, h, f$ is a standard basis of the underlying $\gs\gl_2$. Assume that $\lambda \neq 0$ and consider the irreducible highest weight
$\gG$--module with highest weight vector $v$ such that $e_n \cdot v = 0$ for $n \geq 0$, $f_n \cdot v = 0$ for $n >0$, $h_0 \cdot v = \lambda v$,
$K \cdot v = 0$. Fix $n$ and set $v_k := (e_k f_{-n + k}) \cdot  v$ and $X_k := h_k h_{n-k}$ for $1 \leq k \leq n-1$. A direct computation shows that
$$
X_k \cdot v_l = \left\{ \begin{array} {rcccc}
4 \lambda v & {\text { if }} & l>k & {\text { and }} & n < k+l \\
-4 \lambda v & {\text { if }} & l \leq k & {\text { and }} & n \geq k+l \\
0 & &  &  & {\text { otherwise,}} 
\end{array}
\right.
$$
Taking $1 \leq k, l < \frac{n}{2}$ we see that the corresponding matrix is invertible and hence the dimension of the weight space
with weight $(\lambda, 0) - n \delta$ is at least $\frac{n}{2}$.
\hfill $\square$

\medskip

\subsection{Cuspidal modules} \label{subsec_cuspidal}
Assume that $M \in \cW_\fin$ is a cuspidal module. The result of Britten and Lemire, \cite{BL}, implies that $K$ acts trivially on $M$. 
Theorem \ref{reduce} shows that $M$ is the twisted localization of a parabolically induced $\gG$--module with bounded weight multiplicities.
Propositions \ref{prop4.1} and \ref{prop4.2} show that the corresponding parabolic subalgebra should be imaginary.
Applying Theorem \ref{prop4.1} we see that $M$ is the result of a sequence of twisted localizations applied to a simple weight
module which is parabolically induced from a parabolic subalgebra with Levi component isomorphic to the Heisenberg algebra $\cH = \cL(\gh) \oplus \F D
\oplus \F K$. The irreducible weight modules of $\cH$ (up to a shift of the action of $D$) 
are in a bijection with graded homomorphisms $\Lambda: \gh \otimes \F[t, t^{-1}] \to \F[t, t^{-1}]$. A homomorphism $\Lambda$ is called an exp--polynomial
function if there exists distinct nonzero scalars $\lambda_1, \ldots, \lambda_k \in \F$ and polynomials $p_1, \ldots, p_k$ such that
$$
\Lambda(h \otimes t^n) = (\sum_{i=1}^k p_i(n) \lambda_i^h) t^n.
$$
If $\gP$ is a parabolic subalgebra of $\gG$ with Levi component $\cH$ and $\Lambda: \gh \otimes \F[t, t^{-1}] \to \F[t, t^{-1}]$ is a graded 
homomorphism. We denote by $V_\gP(\Lambda)$ the irreducible $\gG$--module corresponding to $\Lambda$. The following statement is not
difficult to prove.

\medskip

\begin{proposition} \label{prop4.3}
$\phantom{x}$

\noindent
{\rm (i)} $V_{\gP} (\Lambda) \in \cW_\fin$ if and only if $\Lambda$ is an exp--polynomial function.

\noindent
{\rm (ii)} If the dimensions of the weight spaces of $V_\gP(\Lambda)$ are bounded then $\Lambda$ is purely exponential, i.e. all polynomials $p_i$
are constants. In this case $V_\gP(\Lambda)$ is a loop module, i.e. $V_\gP(\Lambda) \cong V_{\lambda_1, \ldots, \lambda_k}(V_1 \otimes \ldots \otimes 
V_k)$ or $V_\gP(\Lambda) \cong V_{\lambda_1, \ldots, \lambda_k}^\sigma(V_1 \otimes \ldots \otimes 
V_k)$ if $\gG$ is twisted, where $V_1, \ldots, V_k$ are irreducible highest weight $\gg$--modules.
Conversely, the weight multiplicities of $V_{\lambda_1, \ldots, \lambda_k}(V_1 \otimes \ldots \otimes 
V_k)$ or $V_\gP(\Lambda) \cong V_{\lambda_1, \ldots, \lambda_k}^\sigma(V_1 \otimes \ldots \otimes 
V_k)$ are bounded if and only if at most one of the modules $V_1, \ldots, V_k$ is infinite dimensional.

\noindent
{\rm (iii)} If $\gG$ admits a simple weight module with bounded weight multiplicties which is not integrable, then $\gG \cong A_l^{(1)}$ or
$\gG \cong C_l^{(1)}$. 
\end{proposition}

\medskip

\noindent
{\bf Sketch of Proof.} Part (i) is proved by a straightforward computation similar to the one in the proof of Proposition \ref{prop4.2} above.
Part (ii) follows from a result of B. Wilson, \cite{W} and the explicit for of the loop modules. Part (iii) follows from (ii). 
\hfill $\square$

\medskip

Combining Theorem \ref{reduce}, Proposition \ref{loc_tensor}, and Proposition \ref{prop4.3} we obtain the main theorem about cuspidal $\gG$ modules.

\begin{theorem} \label{theorem_cuspidal} $\phantom{x}$

\noindent
{\rm (i)} If $\gG$ admits cuspidal modules then $\gG \cong A_l^{(1)}$ or $\gG \cong C_l^{(1)}$. 

\noindent
{\rm (ii)} $M \in \cW_\fin$ is a cuspidal $\gG$--module if and only if $M$ is isomorphic (up to a shift) to 
$$V_{a_0, a_1, \ldots, a_k}(N \otimes V_1 \otimes \ldots \otimes V_k),$$ where
$a_0, \ldots, a_k$ are distinct nonzero scalars, $N$ is a cuspidal $\gg$--module, and $V_1, \ldots, V_k$ are finite dimensional $\gg$--modules.

\noindent
{\rm (iii)} The isomorphism (up to a shift) of cuspidal modules is given by 
$V_{a_0, a_1, \ldots, a_k}(N, V_1, \ldots, V_k) \cong V_{a_0', a_1', \ldots, a_{k'}'}(N', V_1', \ldots, V_{k'}')$ if and only if
$k = k'$, $N \cong N'$, and (after a possible relabelling) $V_i \cong V_i'$, $a_i' = a a_i$ for some nonzero $a \in \F$.
\end{theorem}

\end{document}